\documentclass{article}
\usepackage[]{geometry}
\usepackage{amssymb,amsthm,amsmath,amsxtra,dsfont}
\usepackage[titletoc,title]{appendix}
\usepackage[nottoc]{tocbibind}   % add bibliography to table of contents
\usepackage{graphicx,tikz,esint,pgfplots}
\usepackage{hyperref} \hypersetup{colorlinks=true, citecolor=blue, linkcolor=black}
\usepackage{mathrsfs}
\numberwithin{equation}{section}

\newcommand{\pb}{\bar{\partial}}
\renewcommand{\bar}{\overline}
\newcommand{\supp}{\operatorname{supp}}

\newcommand{\tr}{\operatorname{Tr}}
\newcommand{\Var}{\operatorname{Var}}
\newcommand{\E}{\mathbb{E}}

\newcommand{\Ginibre}{$\infty$-Ginibre }

\newcommand{\1}{\mathbf{1}}
\newcommand{\A}{\mathrm{A}}
\newcommand{\C}{\mathbb{C}}
\newcommand{\Cu}{\operatorname{C}}

\newcommand{\p}{\partial}
\newcommand{\D}{\mathbb{D}}
\renewcommand{\d}{\mathrm{d}}

\newcommand{\h}{\mathfrak{h}}

\newcommand{\J}{\mathbf{J}}
\renewcommand{\k}{\mathbf{k}}
\newcommand{\K}{\mathcal{K}}
\renewcommand{\L}{\mathbf{L}}
\newcommand{\M}{\mathrm{M}}
\newcommand{\N}{\mathbb{N}}
\newcommand{\No}{\mathcal{N}}
\renewcommand{\O}{\underset{N\to\infty}{O}} 
\newcommand{\R}{\mathbb{R}}

\renewcommand{\S}{\mathscr{S}} 

\renewcommand{\u}{\mathrm{u}}
\newcommand{\w}{\mathrm{w}}
\newcommand{\X}{\mathfrak{X}}
\newcommand{\Y}{\mathrm{Y}}
\newcommand{\Z}{\mathrm{Z}}
\newcommand{\z}{\mathrm{z}}

\newtheorem{theorem}{Theorem}[section]

\newtheorem{proposition}[theorem]{Proposition}

\newtheorem{lemma}[theorem]{Lemma}

\theoremstyle{definition} \newtheorem{remark}{Remark}[section]

\title{Incomplete determinantal processes: from random matrix to Poisson statistics}
\date{}
\author{Gaultier Lambert\thanks{Department of Mathematics, University of Zurich, gaultier.lambert@math.uzh.ch}}

\begin{document}

\maketitle

\begin{abstract}
\noindent \normalsize
We study linear statistics of a class of determinantal processes which interpolate between Poisson and GUE/Ginibre statistics in dimension 1 or~2. These processes are obtained by performing an independent Bernoulli percolation on the particle configuration of a log-gas confined in a general potential. 
We show that, depending on the expected number of deleted particles, there is a universal  transition for mesoscopic linear statistics. Namely, at small scales, the point process behave according to random matrix theory, while, at large scales, it needs to be renormalized because the variance of any linear statistic diverges. The crossover is explicitly characterized as the superposition of a $H^{1}$-- or $H^{1/2}$-- correlated  Gaussian noise depending on the dimension and an independent Poisson process. 
The proof consists in computing the limits of the cumulants of linear statistics using the asymptotics of the correlation kernel of the process.
\end{abstract}

  \section{Introduction and results}
 \subsection{Introduction}

In these notes, we consider a log-gas,  also known as $\beta$-ensemble or one component plasma, in dimension $1$ or $2$ at inverse temperature $\beta=2$. Let $\X = \R$  equipped with the Lebesgue measure  $d\mu = dx$ or $\C$ equipped with the  area measure $d\mu= d\A = \frac{r dr d\theta}{\pi}$. Let also $V\in C^2(\X)$ be a real-valued function such that for a $\nu >0$, 
\begin{equation} \label{potential}
 V(z) \ge (1+\nu) \log|z|   
\hspace{.5cm}\text{as}\hspace{.3cm} |z| \to\infty .
\end{equation}
We consider the probability measure on $\X^N$ with  density $G_N(x)=e^{-\beta \mathscr{H}^N_V(x)}/Z^N_V$ where the Hamiltonian is
 \begin{equation} \label{log_gas}
\mathscr{H}^N_V(x)=  \hspace{-.15cm} \sum_{1\le i<j \le N}   \hspace{-.15cm}\log| x_i -x_j|^{-1}  + N \sum_{j=1}^N V(x_j) .
 \end{equation}   
Regardless of the dimension $\d$, the condition $\beta=2$ implies that if a configuration $(\lambda_1,\dots, \lambda_N)$ is sampled from $G_N$, then the point process $\Xi := \sum_{k=1}^N \delta_{\lambda_k}$ is  determinantal  with a correlation kernel
  \begin{equation} \label{kernel}
K^N_V(z,w) = \sum_{k=0}^{N-1} \varphi_k(z)\overline{\varphi_k(w)} , 
\end{equation} 
with respect to $\mu$. Moreover, for all $k\ge 0$, 
\begin{equation} \label{OP}
 \varphi_k(x) = P_k(x) e^{-N V(x)} ,
\end{equation}
 where  $\{P_k \}_{k=0}^\infty$ is the sequence of orthonormal polynomials\footnote{For any $k \ge 0$, $P_k$ is a polynomial of degree $k$ and its leading coefficient is positive.} with respect to the weight $e^{- 2N V(x)}$ on $L^2(\mu)$. 
It turns out that for $\beta=2$,  the density $G_N$ also corresponds to the joint law of the eigenvalues of the ensemble of Hermitian (or normal) matrices with weight $e^{-2N \tr V(M)}$ on $\X=\R$ (or $\X=\C$).   
 In particular, when $V(z)=|z|^2$, these correspond to the well-know Gaussian Unitary (GUE) and Ginibre ensembles respectively. 
It is well known that if the condition \eqref{potential} holds, the thermodynamical limit of the log-gas is described by an equilibrium measure which has compact support. 
Moreover, if the potential $V\in C^2(\X)$,  then the equilibrium measure is absolutely continuous and we let $\varrho_V$ be its density. This implies that for any bounded test function $f \in C(\X)$, as $N\to+\infty$, 
\begin{equation} \label{LLN}
\frac{1}{N}\E\big[\Xi(f)\big]  = \int f(x) u^N_V(x) d\mu(x) \to \int f(x) \varrho_V(x) d\mu(x) ,
\end{equation}
where the expected density of states is given by $u^N_V(x) = N^{-1}K^N_V(x,x)$. 
The asymptotics \eqref{LLN} follows either from potential theory for general $\beta>0$ or from the asymptotics of the correlation kernel \eqref{kernel} when $\beta=2$. 
We refer to \cite[Section~2.6]{AGZ} for a proof of the large deviation principle in dimension~1 and to \cite{CHM16} for analogous results  for Coulomb gases in higher dimension and further references.

\medskip 
 
%For large classes of test functions, the asymptotic properties of linear statistics of these ensembles have been investigated and laws of large numbers (LLN) and central limit theorems (CLT) have been proved in various regime, see \cite{Johansson98, AHM11, AHM15,  BD14, BD16, Lambert_b, BL16} to cite of few of the main contributions or the books \cite{Forrester10, PS11} for a review of the theory. \\

 In the following, we  consider the problem of describing the fluctuations of the so-called \emph{thinned}  log-gases in dimension $\d =1,2$.   
In general, a \emph{thinned} or \emph{incomplete  point  process} is defined by performing a Bernoulli percolation on the configuration of a the original process. 
That is  the incomplete log-gas, denoted by $\widehat\Xi$, is obtained by deleting independently each particle with probability $q_N \in (0,1)$ or by keeping it  with probability $p_N=1-q_N$. 
It turns out that the incomplete process  $\widehat\Xi$ is also determinantal with correlation kernel 
$\widehat{K}^N_V(z,w) = p_NK^N_V(z,w)$; see the appendix~\ref{A:kernel} for a short proof.   In the context of random matrix theory,  this procedure was first considered by Bohigas and Pato \cite{BP04, BP06} who showed that it gives rise to a crossover to Poisson statistics  and the problem of  rigorously analyzing this transition in the context of Coulomb gases was popularized by Deift in \cite[Problem 2]{Deift17}. 
Indeed, these types of transitions are supposed to arise in many different contexts in statistical physics, such as the localization/delocalization phenomena, the crossover from chaotic to regular dynamics in the theory of quantum billiards, or in the spectrum of certain band random matrices, see \cite{Spencer11, EK14a, EK14b} and reference therein.
Although such transitions are believed to be non-universal, the model of Bohigas and Pato is arguably one of the most tractable  to study this phenomenon because it is determinantal.
In a different context, the effect of thinning determinantal process on statistical inferences has been recently discussed in~\cite{LMR_15} and it should be emphasized that the general strategy explained in section~\ref{sect:cumulant} applies to more general determinantal processes, see theorem~\ref{thm:main}.  For instance, our method  applies to the Sine and the \Ginibre processes which describes the local limits of the log-gases in dimensions $1$ and $2$ respectively. In fact,  this paper is motivated by an analogous result obtained recently by Berggren and Duits for smooth linear statistics of the incomplete Sine and CUE processes \cite{Berggren17}. 
Based on the fact that these processes come from integrable operators,  they fully characterized the transition for a large class of mesoscopic linear statistics and suggested that it should be universal for thinned point processes coming from random matrix theory. 
There are also results for the gap probabilities of the critical thinned ensembles. In~\cite{BDIK15, BDIK16}, for the Sine process, Deift et al.~computed very detailed asymptotics for the crossover from the {\it Wigner surmise} to the exponential  distribution  making rigorous a prediction of Dyson \cite{Dyson95}, and Charlier--Claeys obtained an analogous result for the CUE \cite{CC17}. 
The contribution of this paper is to elaborate on universality
for smooth linear statistics of  $\beta$-ensemble in dimension $1$ or $2$ when $\beta=2$.  Although our proof relies on the determinantal structure of these models, instead of  the connection with Riemann-Hilbert problems used in the previous works, we apply the cumulants method which appears to be very robust to study the asymptotic fluctuations of smooth linear statistics.  

\medskip

Let us point out that based on the theory of \cite{HKPV06}, an alternative correlation kernel for the incomplete process is 
\begin{equation} \label{kernel_2}
\widehat{K}^N_V(z,w) = \sum_{k=0}^\infty J_k^N  \varphi_k(z) \varphi_k(w) , 
\end{equation}
where  $(J_k^N)_{k=1}^\infty$ is a sequence of i.i.d.~Bernoulli random variables with expected values $\E[J_k^N]= p_N \1_{k<N} $.  
This shows that removing particles builds up randomness in the system and when the disorder becomes sufficiently strong, it will behave like a Poisson process rather than according to random matrix theory. 
%Although, this heuristics is difficult to turn into a proof, especially in the mesoscopic regime.  

\medskip

To keep the analysis as simple as possible, we will restrict  ourself to real-analytic~$V$, although the results should be valid for more general potential as well (especially in dimension~1 where the asymptotics of the correlation kernels have been studied in great generality). We  keep track of the transition by looking at linear statistics $\widehat{\Xi}(f) = \sum f(\lambda)$ for smooth test functions, where the sum is over the configuration  of the incomplete log-gas.   
 The random matrix regime is characterized by the property that the fluctuations of $\widehat{\Xi}(f)$ are of order 1 and described by a {\it universal Gaussian noise} as the number of particles tends to infinity. On the other hand, in the Poisson regime, the variance of any non-trivial statistic diverges and, once properly renormalized, the point process converges in distribution to a {\it white noise}. 
In the remainder of this introduction, we first formulate our assumptions and main results  for the fluctuations of the incomplete 2-dimensional log-gases. Then, we will present  analogous results in dimension~1.     

\medskip

In what follows, we let $\mathcal{C}^k_c(\S)$ be the set of functions with $k$ continuous derivatives and compact support in $\S \subset \X$ and we use the notation:
$$
\p = (\p_x - i \p_y)/2 \ , \hspace{.5cm} \pb = (\p_x + i \p_y)/2 \hspace{.5cm}\text{and}\hspace{.4cm}
\Delta = \p\pb .
$$
If $\varrho_V$ is the equilibrium density function, we also denote 
\begin{equation}\label{bulk}
\S_V:=  \big\{ x\in\X : \varrho_V(x) >0 \big\} . 
\end{equation}

\subsection{Main results for 2-dimensional Coulomb gases} \label{sect:gas2}

If the potential $V$ is real-analytic and satisfies the condition \eqref{potential},
then the log-gas lives on the compact set $\overline{\S_V} \subset \C$ which is called the {\it droplet} and the equilibrium density is given by $\varrho_V = 2\Delta V \1_{\S_V}$. % Note that this density may vanish inside the droplet and we shall focus on describing the fluctuations inside the  {\it bulk} $\S_V$, \eqref{bulk}.
 It is also well known that the bulk fluctuations of a two dimensional log-gas around its equilibrium configuration are described by a centered Gaussian process $\mathrm{X}$ with correlation structure: 
\begin{equation} \label{H1_noise}
\E\big[\mathrm{X}(f) \mathrm{X}(g) \big] = \frac{1}{4}  \int_\C  \nabla f(z) \cdot \nabla g(z)\ d\A(z) 
=  \int_\C  \p f(z) \pb g(z) d\A(z)    , 
\end{equation}
for any (real-valued) smooth functions $f$ and  $g$.
Modulo constants, the RHS~of formula \eqref{H1_noise} defines a Hilbert space, denoted $ H^{1}(\C)$, with norm:
\begin{equation} \label{H1_norm}
\|f \|_{H^1(\C)}^2 =  \int_\C  \left| \p f(z) \right|^2 d\A(z) .    
\end{equation}
Therefore the stochastic process $\mathrm{X}$ is called  a $H^1$-Gaussian noise. 
 The  central limit theorem (CLT)  was first established for the Ginibre process by Rider and Vir\`{a}g \cite{RV07a} for $\mathcal{C}^1$ test functions with at most exponential growth at infinity.
For general real-analytic potentials, it was proved in \cite{AHM15} that for any smooth function $f$ with compact support, one has as $N\to\infty$,
\begin{equation} \label{clt_1}
\Xi(f) - \E\big[\Xi(f) \big]  \Rightarrow  \mathrm{X}(f^\dagger)
\end{equation}
 where $f^\dagger$ is the (unique) continuous and bounded function on $\C$ such that  $f^\dagger= f$ on the droplet $\overline{\S_V}$ and $\Delta f^\dagger =0$ on $\C\backslash\overline{\S_V}$. Actually, when $\supp(f)  \subset \S_V $,  we have  $f^\dagger = f$ on $\C$ and the CLT was obtained previously in the paper \cite{AHM11} from which part of our method is inspired. We also refer to \cite{BBNY16, LS18} for more recent proofs which hold for general $\beta>0$. 
By convention,  in \eqref{H1_noise} and below, $\Rightarrow$ means that the convergence holds in distribution and that all moments of the random variable converge.

 \medskip
  
 In order  to describe the crossover from the  $H^1$-Gaussian noise to white noise, let $\Lambda_\eta$ be a mean-zero Poisson process with intensity $\eta \in L^\infty(\X)$. This process is characterized by the fact that for any function $f\in \mathcal{C}_c(\X)$, the Laplace transform of the random variable $\Lambda_\eta(f)$ is well-defined and given by
\begin{equation}\label{Poisson}
\log \E\big[ \exp \Lambda_{\eta}(f) \big] = \int_{\X} \big( e^{f(z)} - 1 - f(z) \big) \eta(z) d\mu(z) . 
\end{equation}

%If $X$ and $Y$ are two real-valued random variables, we denote by $X \oplus Y$ the law of the $X+Y$. 

\begin{theorem} \label{thm:crossover_2}
% Let $\mathrm{X}(f)=  \No\big(0, \| f \|_{H^{\d/2}(\X)}^2\big)   $. 
Let $\mathrm{X}$ be a $H^1$-Gaussian noise and $\Lambda_{\tau\varrho_V}$ be an independent  Poisson process with intensity $\tau\varrho_V$ where $\tau>0$ defined on the same probability space.
Let  $f \in \mathcal{C}^3_c(\S_V)$, $p_N = 1-q_N$, and let  $T_N = N q_N$. 
As $N\to\infty$ and $q_N \to 0$, we have
\begin{align}
\label{Normal_2}
\widehat{\Xi}(f) - \E\big[\widehat{\Xi}(f) \big]  &\Rightarrow  \mathrm{X}(f) &\text{if }\ T_N \to 0 , \hspace{.2cm}\\ 
\label{Poisson_2}
\frac{\widehat{\Xi}(f) - \E\big[\widehat{\Xi}(f) \big]}{\sqrt{T_N} } &  \Rightarrow
 \No\left(0, \int  f(z)^2 \varrho_V(z) d\A(z) \right)   
&\text{if }\ T_N \to \infty ,\\
\label{crossover_2}
\widehat{\Xi}(f) - \E\big[\widehat{\Xi}(f) \big] & \Rightarrow  \mathrm{X}(f) + \Lambda_{\tau\varrho_V}(-f) &\text{if }\ T_N \to \tau  . \hspace{.25cm}
\end{align}
\end{theorem}

The proof of theorem~\ref{thm:crossover_2}  is based on the cumulants' method and it is explained in details in section~\ref{sect:cumulant}.
In particular, we formulate a result -- theorem~\ref{thm:main} -- valid for general determinantal point process and might be of independent interest.
 The details of the proof of theorem~\ref{thm:crossover_2}  are given in section~\ref{sect:global2}.  
Our method relies on the approximations of the correlation kernel $K^N_V$ from \cite{AHM11}   -- see lemma~\ref{thm:Bergman_kernel} below -- and it restricts us to work with test functions which are supported inside the bulk. However, the result should be true for general functions  if we replace $f$ by $f^\dagger$ on the RHS~of \eqref{Normal_2} and \eqref{crossover_2}.  
%In fact, for the Ginibre ensemble, $V(z)=|z|^2$, or other  rotationally invariant potential, one can also apply the combinatorial method of \cite{RV07a}    
%to prove the counterpart of theorem~\ref{thm:crossover_2} for any test function $f(z)$ which is a bi-variate polynomial in $z$ and $\bar{z}$. 

\medskip

Theorem~\ref{thm:crossover_2} can be interpreted as follows.
In the regime $Nq_N \to 0$, virtually no particles are deleted and linear statistics behave according to random matrix theory. 
On the other hand, in the regime $Nq_N \to \infty$, the variance of a linear statistic diverge. So, if we renormalize the random variable $\widehat{\Xi}(f)$, we obtain a classical CLT and  \eqref{Poisson_2} shows that the limit is described by a white noise supported on $\S_V$ whose intensity is  the equilibrium measure~$\varrho_V$. 
In the critical regime, when the expected number of deleted particles equals $\tau>0$, the limiting process is the superposition of a $H^1$-correlated Gaussian noise and an independent mean-zero Poisson process  applied to $-f$.  
Finally, by using formula \eqref{Poisson}, it is not difficult to check that as $\tau\to\infty$, the random variable  
$$
\frac{\Lambda_{\tau\varrho_V}(-f)}{\sqrt{\tau}}  \Rightarrow  \No\left(0, \int  f(z)^2 \varrho_V(z) d\A(z) \right)   , 
$$
  so that the critical regime clearly interpolates between  \eqref{Normal_2} and~\eqref{Poisson_2}. 
 
\medskip
 
In fact, the crossover is more interesting at mesoscopic scales. 
Namely, the density of a log-gas is of order $N$  and one can also investigate fluctuations at small scales by zooming inside the bulk of the process. If $L_N \nearrow \infty$,  $x_0 \in \S_V$, and $f\in \mathcal{C}_c(\C)$, we consider the test function
\begin{equation}  \label{meso}
f_N(z) = f\big(L_N(z-x_0) \big) .
\end{equation}

The regime $L_N=N^{1/\d}$ is called microscopic and it was shown in~\cite[Proposition~7.5.1]{AHM11} that  when $\d=2$, 
$$\Xi (f_N) \Rightarrow \Xi^\infty_{\varrho_V(x_0)}(f) , $$
where the process $ \Xi^\infty_{\rho}$ is called the {\bf \Ginibre process} with density $\rho>0$. 
It is a determinantal process on $\C$ with correlation kernel
\begin{equation}\label{Ginibre}
K^\infty_\rho(z,w) =  \rho e^{\rho(2z\bar{w} - |z|^2 -|w|^2 )/2} . 
\end{equation}
Based on the argument from \cite{AHM11}, it is straightforward to verify that the incomplete process has a local limit as well:
\begin{equation}\label{local_limit_1}
\widehat{\Xi} (f_N) \Rightarrow \widehat{\Xi}^{\infty}_{\varrho_V(x_0) ; p}(f)
\hspace{.6cm}\text{as}\hspace{.4cm}p_N\to p \hspace{.4cm}\text{and}\hspace{.4cm} N\to\infty. 
\end{equation}
For any $0<p\le 1$, $ \widehat{\Xi}^{\infty}_{\varrho; p}$ is  a (translation invariant) determinantal process on $\C$ with correlation kernel  $p K^\infty_\varrho(z,w) $.  This process  is constructed by running an independent Bernoulli percolation with parameter $p$ on the point configuration of  the \Ginibre process with density $\rho>0$. 
In particular, \eqref{local_limit_1} shows that one needs to delete a non-vanishing fraction of the $N$ particles of the gas in order to get a local limit which is different from random matrix theory. 
It was proved in \cite{RV07b} that, as the density $\rho\to\infty$, the fluctuations of the \Ginibre process are of order 1 and described by the  $H^1$-Gaussian noise:
 $$
 \Xi^\infty_\rho(f) - \rho \int_{\C} f(z) d\A(z)\  \Rightarrow \mathrm{X}(f) 
 $$
for any $f \in H^1\cap L^1(\C)$. 
Therefore it is expected that, in the mesoscopic regime, $\text{i.e. }L_N = o(\sqrt{N})$, the asymptotic fluctuations of the linear statistic  $\Xi(f_N)$ are universal and described by  $\mathrm{X}(f)$.  
However, to our best knowledge, a proof was missing from the literature and, in section~\ref{sect:2}, we show that this fact follows quite simply  by combining the ideas from  \cite{RV07b} and \cite{AHM11}.

\begin{theorem} \label{meso:2}
 Let  $x_0 \in \S_V$, $f \in \mathcal{C}_c^3(\C)$, $\alpha \in (0,1/2)$, and let $f_N$ be given by formula \eqref{meso} with  $L_N =N^\alpha$. Then, we have  as $N \to\infty$,
  $$
 \Xi(f_N) - \E\big[ \Xi(f_N)\big] \Rightarrow \mathrm{X}(f) .
 $$
 \end{theorem}

 Using the same method, we can also describe the fluctuations of smooth mesoscopic linear statistics of a incomplete  Coulomb gas. 

\begin{theorem} \label{thm:crossover_3}
% Let $\mathrm{X}(f)=  \No\big(0, \| f \|_{H^{\d/2}(\X)}^2\big)   $.
Let $\mathrm{X}$ be a  $H^1$-Gaussian noise and $\Lambda_{\tau}$ be an independent Poisson process with constant intensity $\tau>0$ on $\C$.
Let  $x_0 \in \S_V$,  $f \in \mathcal{C}^3_c(\S_V)$, $\alpha \in (0,1/2)$, and let $f_N$ be the mesoscopic test function given by formula \eqref{meso} with  $L_N =N^\alpha$. We also let  $p_N = 1-q_N$ and $T_N= N q_N  L_N^{-2} \varrho_V(x_0)$. We have  as $N\to\infty$ and $q_N\to0$, 
\begin{align*}
%\label{Normal_3}
\widehat{\Xi}(f_N) - \E\big[\widehat{\Xi}(f_N) \big]  &\Rightarrow  \mathrm{X}(f) &\text{if }\ T_N \to 0 , \hspace{.2cm}\\ 
%\label{Poisson_3}
\frac{\widehat{\Xi}(f_N) - \E\big[\widehat{\Xi}(f_N) \big]}{\sqrt{T_N} } &  \Rightarrow
 \No\left(0, \int  f(z)^2 d\A(z) \right)   
&\text{if }\ T_N \to \infty ,\\
%\label{crossover_3}
\widehat{\Xi}(f_N) - \E\big[\widehat{\Xi}(f_N) \big] & \Rightarrow  \mathrm{X}(f) + \Lambda_{\tau}(-f) &\text{if }\ T_N \to \tau  .\hspace{.25cm}
\end{align*}
%as $N\to\infty$, where $\tilde{t} =  \varrho_V(x_0) t$. 
\end{theorem}

The proof of theorem~\ref{thm:crossover_3}  follows the same strategy are that of theorem~\ref{thm:crossover_2} and the technical differences are explained in section~\ref{sect:global2}.
This result shows that, at mesoscopic scales, the transition occurs when the  mesoscopic density of deleted particles which is given by the parameter $T_N>0$ converges to a positive constant $\tau$.
In contrast to previous results, this transition appears to be non-Gaussian and it is somewhat surprising that it can also be described in an elementary way. In dimension 1,  one can obtain a crossover from GUE eigenvalues to a Poisson process  by letting  independent points evolved according to  Dyson's Brownian motion. This leads to a determinantal process sometimes called the {\it deformed GUE} whose kernel depends on the diffusion time, see \cite{Johansson01}. 
This model also exhibit a  transition which has been analyzed for mesoscopic linear statistics in  \cite{DJ16} and  it was proved that 
 the critical fluctuations are Gaussian. One can also consider non-intersecting Brownian motions on a cylinder. It turns out that this point process describes the positions of free fermions confined in a harmonic trap at a temperature $\tau>0$. 
 It was established in \cite{Johansson07}, see also \cite{DDMS14}, that the corresponding grand canonical ensemble is determinantal with a correlation kernel of the form \eqref{kernel_2} with  for $k\ge 0$, 
 $$\E[J_k^N] = \frac{1}{1+ \exp(\frac{k-N}{\tau})} . $$  
 For sufficiently small temperature, this system behaves like its ground-state, the GUE, while it
behaves like a Poisson process at larger temperature.
 It was proved in \cite{JL15} that this leads to yet another crossover where non-Gaussian fluctuations are observed at the critical temperature. However,  to the author's knowledge, in contrast to the incomplete ensembles considered here, the critical processes discovered  in \cite{JL15}  cannot be described in simple terms like theorem~\ref{thm:crossover_4} below.
 %%% Work

 \subsection{Main results for eigenvalues of unitary invariant Hermitian random matrices} \label{sect:gas1}
 
For 1-dimensional log-gases, for general $\beta>0$ and for a large class of potentials, Johansson established in \cite{Johansson98} the existence of the equilibrium measure  and also managed to describe the fluctuations around the equilibrium configuration.
To state the result in a universal way, note that one can make an affine rescaling of the potential  and assume that  $\mathscr{I}_V \subset [-1,1] $. 
 If $V$ is a polynomial and $\mathscr{I}_V =(-1,1)$,  Johansson proved that linear statistics of the  process $\Xi$ 
 satisfy a central limit theorem:
\begin{equation} \label{Johansson}
\Xi(f) - N \int_\R f(x) \varrho_V(x) dx \Rightarrow \Y(f) \hspace{.5cm}\text{as }N\to\infty ,
\end{equation}
for any $f \in \mathcal{C}^2(\R)$ such that $f'(x)$ grows at most polynomially as $|x|\to\infty$. 
The process $\Y$ is a centered Gaussian noise defined on $[-1,1]$ with covariance structure:
\begin{equation} \label{Chebychev_covariance}
\E\big[ \Y(f)\Y(g) \big] = 
 \frac{1}{4}  \sum_{k=1}^\infty k \mathrm{c}_k(f)  \mathrm{c}_k(g) . 
 \end{equation}
In \eqref{Chebychev_covariance}, $\mathrm{c}_k(f)$ denote the Fourier-Chebyshev coefficients of the function $f$:
\begin{equation}\label{Fourier_T}
\mathrm{c}_k(f)=\frac{2}{\pi} \int_{-1}^1 f(x) T_k(x) \frac{dx}{\sqrt{1-x^2}} , 
\end{equation}
 where $(T_k)_{k=0}^\infty$ are  the Chebyshev polynomials of the first kind\footnote{$T_k(\cos \theta) = \cos(k\theta)$ for any $k\ge 0$ and $\theta\in\R$.}. 
The CLT \eqref{Johansson}  holds for more general potentials and for other  orthogonal polynomial ensembles as well, see \cite[section~11.3]{PS11} or \cite{BG13a, BD17, Lambert_b} and   
it is  known that {\it the one-cut condition}, i.e.~the assumption that the support of the equilibrium measure is connected, is necessary. Otherwise, the asymptotic fluctuations of a generic linear statistic $\Xi(f)$ are still of  order~1 but are not Gaussian, see \cite{Pastur06, Shcherbina13, BG13b}. In fact, the one-cut condition is closely related to the fact  the recurrence  coefficients (see formula \eqref{recurrence} below) which defines the orthogonal polynomials $(P_k)_{k\ge0}$ appearing in the correlation kernel \eqref{kernel} satisfy for any $j\in\mathbb{Z}$,  
\begin{equation}\label{recurrence_limit} \lim_{N\to\infty} a^N_{N+j} = 1/2 \hspace{.5cm}\text{and}\hspace{.5cm}
   \lim_{N\to\infty}b^N_{N+j} = 0;
   \end{equation}
   see remark~\ref{rk:one_cut} below.
  Like for 2-dimensional Coulomb gas, we obtain analogous transitions for the eigenvalues of random unitary invariant Hermitian matrices.

\begin{theorem} \label{thm:crossover_1}
 Let $p_N = 1-q_N$, $T_N = N q_N$, and suppose that the recurrence coefficients of the orthogonal polynomials $\{P_k \}_{k=0}^\infty$  satisfy the conditions \eqref{recurrence_limit}.  
 Then, for any polynomial $Q$,  we obtain as $N\to\infty$ and $q_N \to 0$, 
\begin{align*}
%\label{Normal_1}
\widehat{\Xi}(Q) - \E\big[\widehat{\Xi}(Q) \big]  &\Rightarrow  \mathrm{Y}(Q) &\text{if }\ T_N \to 0 , \hspace{.2cm}\\ 
%\label{Qoisson_1}
\frac{\widehat{\Xi}(Q) - \E\big[\widehat{\Xi}(Q) \big]}{\sqrt{T_N} } &  \Rightarrow
 \No\left(0, \int_{\R}  Q(x)^2 \varrho_V(x) dx \right)   
&\text{if }\ T_N \to \infty ,\\
%\label{crossover_1}
\widehat{\Xi}(Q) - \E\big[\widehat{\Xi}(Q) \big] & \Rightarrow  \mathrm{Y}(Q) + \Lambda_{\tau\varrho _V}(-Q) &\text{if }\ T_N \to \tau  , \hspace{.25cm}
\end{align*}
where the Poisson process $ \Lambda_{\tau\varrho_V}$ is independent from the Gaussian process $\Y$ and both are defined on $\mathscr{I}_V =(-1,1)$. 
\end{theorem}

The proof of theorem~\ref{thm:crossover_1} is also based on the cumulants' method and on theorem~\ref{thm:main}. However, the technical details, which are explained in sections~\ref{sect:aop} and~\ref{sect:glob1}, rely on the formulation from \cite{Lambert_b} and are very different from that of the proof of theorem~\ref{thm:crossover_2}.
%Let us also mention that in the macroscopic regime, we give a combinatorial proof which is based on the asymptotics of the recurrence coefficients of the orthogonal polynomials \eqref{OP} and can be generalized to the multi--cut situation as well as to other  biorthogonal ensembles; see remark~\ref{rk:biorthogonal}. 

%In particular, theorem~\ref{thm:crossover_1} is restricted to the class of polynomial test functions since its proof which is inspired from \cite{Lambert_b} and relies on the formulation of the cumulants in terms of recurrence coefficients; see section~\ref{sect:1} for further details. 
%However, by using certain variance estimates and density arguments in the spirit of \cite[Lemma~2.1]{SW_13}, it should not be difficult to extend the results to a larger class of test functions. 
 %
 
 \medskip
 
 We also obtain the counterpart of theorem~\ref{thm:crossover_1} for mesoscopic linear statistics.
For any function $f\in L^1(\R)$, we define its Fourier transform :
$$
\hat{f}(u) = \int_\R f(x) e^{-2\pi i x u} dx .
$$
We let  $\Z$ be a mean-zero Gaussian process on $\R$ with correlation structure:
\begin{equation}
\E\big[ \Z(h)\Z(g) \big] = \int_{0}^\infty u \hat{h}(u) \overline{\hat{g}(u)} du . 
\end{equation} 
Since $\Var\big[ \Z(f)\big] = \| f\|_{H^{1/2}(\R)}$, the process $\Z$ is usually called the $H^{1/2}$--Gaussian noise. It describes the mesoscopic fluctuations of the eigenvalues of Hermitian random matrices, see \cite{BD16, Lambert_a, HK16}, as well as  the mesoscopic fluctuations of the log-gases for general $\beta>0$, \cite{BL16}, and of certain random band matrices in the appropriate regime \cite{EK14a, EK14b}.

 \begin{theorem} \label{thm:crossover_4}
We let $x_0\in\mathscr{I}_V$, $f\in \mathcal{C}^2_c(\R)$, $\alpha\in (0,1)$,  and $f_N(x) = f\big(N^\alpha(x-x_0) \big)$. We also let $p_N = 1-q_N$ and  $T_N =  q_NN^{1-\alpha} \varrho_V(x_0)$.
 We obtain  as $N\to\infty$ and $q_N\to0$, 
\begin{align}
\label{Normal_4}
\widehat{\Xi}(f_N) - \E\big[\widehat{\Xi}(f_N) \big]  &\Rightarrow \Z(f) &\text{if }\ T_N \to 0 , \hspace{.2cm}\\ 
\label{f_Noisson_4}
\frac{\widehat{\Xi}(f_N) - \E\big[\widehat{\Xi}(f_N) \big]}{\sqrt{T_N} } &  \Rightarrow
 \No\left(0, \int_{\R}  f(x)^2 dx \right)   
&\text{if }\ T_N \to \infty ,\\
\label{crossover_4}
\widehat{\Xi}(f_N) - \E\big[\widehat{\Xi}(f_N) \big] & \Rightarrow \Z(f) + \Lambda_{\tau}(-f) &\text{if }\ T_N \to \tau  , \hspace{.25cm}
\end{align}
where the Poisson process $\Lambda_{\tau}$ has constant intensity $\tau>0$ on $\R$ and  is independent from the  $H^{1/2}$--Gaussian noise $\Z$.  
\end{theorem} 
 
 The proof of theorem~\ref{thm:crossover_4} is quite similar to that of theorem~\ref{thm:crossover_3}. It follows the strategy explained in section~\ref{sect:cumulant} and is based on  the asymptotics for the correlation kernel $K^N_V$ in terms of the sine-kernel, see~\cite{Kuijlaars11}. 
 The details  relies on the method from~\cite{Lambert_a} and are given in section~\ref{sect:meso1}.

 % and the method of~\cite{Lambert_a}. This reduces the proof  to carefully apply the classical Riemann--Hilbert asymptotics  of \cite{Deift99+}. Note that the steepest descent method of \cite{Deift99+} has also been established for many other invariant ensembles, see~\cite{Kuijlaars11}, for which it is also possible to extract the asymptotics  \eqref{sine_1} for the correlation and deduce the analogue of theorem~\ref{thm:crossover_4}.

\subsection{Overview of the rest of the paper}

In section~\ref{sect:cumulant}, we present the strategy of the proofs of the results from sections~\ref{sect:gas2} and \ref{sect:gas1}. We begin by reviewing Soshnikov's cumulants' method.
Then, we explain how to apply it to the incomplete ensemble $\widehat{\Xi}$ and we obtain a general result -- theorem~\ref{thm:main} -- which characterizes the transition from Gaussian to Poisson statistics for general determinantal point processes. 
The rest of the paper consists in verifying the assumptions of theorem~\ref{thm:main} for determinantal log-gases in dimensions 1 and 2. 
In section~\ref{sect:2}, we prove theorems~\ref{thm:crossover_2}--\ref{thm:crossover_3}  for the 2-dimensional log-gases. 
The proof relies on estimates for the correlation kernel \eqref{kernel} which come from the paper~\cite{AHM11} and are collected in the appendix~\ref{A:damping}. 
In section~\ref{sect:1}, we provide the details of the  proofs of theorem~\ref{thm:crossover_1} and~\ref{thm:crossover_4} by relying on the method from \cite{Lambert_b} and \cite{Lambert_a} respectively.

\medskip

In the following,  $C>0$ denotes a numerical constant which changes from line to line.  For any $n\in\N$, we let for all $\mathrm{x} \in\X^n$, 
$$
 d\mu^n(\mathrm{x}) = d\mu(x_1) \cdots d\mu(x_n) . 
$$
   If $(\mathrm{u}_N)$ and $(\mathrm{v}_N)$ are two sequences, we use the notation:  $$\mathrm{u}_N \simeq \mathrm{v}_N 
  \hspace{.5cm}\text{if} \hspace{.5cm} \lim_{N\to\infty}( \mathrm{u}_N - \mathrm{v}_N) =0  .
  $$

%Since the process has a correlation kernel $p_NK^N_V(x,y)$, for any nice test function,

  \subsection*{Acknowledgements}
 I thank Tomas Berggren and Maurice Duits from which I learned about the model of Bohigas and Pato,   for sharing their inspiring work with me and for the valuable discussions which followed. I also thank Mariusz Hynek for many interesting discussions, as well as the referee whose comments help to improve the structure of this paper.
  
  \section{Outline of the proof} \label{sect:cumulant}

In this section, we consider a general state space $\X$ which is a complete separable metric space equipped with a Radon measure $\mu$ as in \cite{Soshnikov, HKPV06} and let $\Xi$ be a sequence of determinantal point processes on $\X$ with correlation kernels $K^N$ which are  reproducing:
\begin{equation} \label{projection}
\int_\X K^N(z, x) K^N(x, w) d\mu(x) = K^N(z,w) . 
\end{equation}
One may think of the parameter $N \in \N$ as the \emph{density} of particles.
 Since it is generally  the case in the context of random matrix theory, we shall also assume that the kernels  $K^N$ are continuous  on $\X\times \X$, Hermitian symmetric,  and that they define locally trace-class integral operators acting on $L^2(\X,\mu)$. 
 The cumulants' method  to analyze the asymptotic distribution of linear statistic of determinantal processes goes back to the work of  Costin and Lebowitz  \cite{CL95}  for count statistics of  the Sine process. The general theory was developed by Soshnikov in \cite{Soshnikov, Soshnikov00, Soshnikov01}  and subsequently applied to many different ensembles coming from random matrix theory, see for instance \cite{RV07b, RV07a, AHM11, BD16, BD17, JL15, Lambert_a, Lambert_b}.
%
  %  
%Then this method has been applied to analyze the  asymptotic fluctuations of various ensembles coming from random matrix theory or other probabilistic models, %
%
 In this section, we show how to implement the cumulants' method to describe the asymptotics law of linear statistics of the incomplete ensemble~$\widehat{\Xi}$ with correlation kernel $p_N K^N(z,w)$ when  the density of particles $0<p_N<1$ converges to 1 in the large $N$ limit. 

   \medskip
   
Let 
  $$
  \mho= \bigcup_{l=1}^\infty \big\{ \k=(k_1,\dots, k_l) \in \N^l  \big\} 
  $$
  and let $\ell(\k)= l$ denote the length of the tuple $\k$. Let us denote the set of compositions of the integer $n>0$ by
  $$ \big\{ \k \vdash n \big\} = \big\{ \k\in\mho :  k_1+ \cdots + k_l = n \} . $$ 
    %$ if $|\k| = k_1+ \cdots + k_l = n $. 
  % In particular,  $\big\{\k \in   \mho : \k \vdash n   \big\}$ is the set of compositions of the integer $n\in\N$.
We also denote by $n\in \mho$ the trivial composition. 
  For any map $\Upsilon : \mho \mapsto \R $, for any function $f: \X \to \R$, and for any $n\in\N$, we define   for all $x\in \X^n$,
\begin{equation}\label{Upsilon_f}
  \Upsilon^n[f](x) = \sum_{\k \vdash n} \Upsilon(\k)\hspace{-.2cm} \prod_{1 \le j\le \ell(\k)} \hspace{-.2cm}f(x_j)^{k_j} . 
 \end{equation}
If $\k \vdash n$, we let $\displaystyle \M(\k) =\frac{n!}{k_1!\cdots k_l!}$ be the corresponding multinomial coefficient and for all integers $n\ge 1$ and $ m \in\{0, \dots , n\}$, we define the coefficients 
\begin{equation}\label{Gamma}
  \gamma^n_m = \sum_{ \k \vdash n}\frac{(-1)^{\ell(\k)}}{\ell(\k)} {\ell(\k) \choose m} \M(\k) .  
\end{equation}
We will also use the notation:
$\displaystyle\delta_{k}(n) = \begin{cases} 1 &\text{if}\ n=k \\ 0 &\text{else} \end{cases}$ for any $k\in\mathbb{Z}$.

 \begin{lemma} \label{thm:combinatorics}
  For all $n\in \N$, we have $\gamma^n_0 = \delta_{1}(n)$ and $  \gamma^n_1 = (-1)^n$.
  \end{lemma}
\proof  The coefficients \eqref{Gamma}  have the generating function:
  $$
  \sum_{n=1}^\infty \sum_{m=0}^n   \gamma^n_m  \frac{x^n q^m}{n!}  = - \log\big(1+(1+ q)(e^x-1)\big) . 
  $$
  In particular, setting $q=0$, we see that $  \gamma^n_0 = \delta_{1}(n) $. 
Moreover, since
 $$
 \sum_{n=1}^\infty   \gamma^n_1  \frac{x^n}{n!}    
 =\left.  \frac{d}{dq}\log\big(1+(1+ q)(e^x-1)\big) \right|_{q=0} =1-e^{-x} , 
 $$
  we also see that $  \gamma^n_1 = (-1)^n$ .\qed\\

Given a test function  $f: \X \to \R$, say locally integrable with compact support, the cumulant generating function of the random variable $\Xi(f)$ is 
$$
\log \E\big[ \exp ( \lambda \Xi(f)) \big] = \sum_{n=1}^\infty \frac{\lambda^n}{n!} \Cu^n_{K^N}[f] .
$$
It was proved by Soshnikov that under our general assumptions,  the cumulants $\Cu^n_{K}[f] $ characterize the law of the linear statistics  $\Xi(f)$ and  that for any $n\in\N$, 
\begin{equation} \label{Soshnikov}
\Cu^n_{K^N}[f]
= - \sum_{l=1}^n   \frac{(-1)^l}{l} \sum_{\begin{subarray}{c} \k \vdash n \\ \ell(\k) =l \end{subarray}} \M(\k) 
\underset{x_{0} =x_l}{\int_{\X^l}} f(x_1)^{k_1} \cdots  f(x_l)^{k_l} \prod_{1\le j\le l} K^N(x_j, x_{j-1}) d\mu^l(x) . 
\end{equation}
Under stronger assumptions, for instance if the kernel $K^N$ has finite rank, this formula makes sense also for 
test functions which are not necessarily compactly supported. 
We use the convention that the variables $x_0$ and $x_l$ are identified in the previous integral. Since we assume that the correlation kernel $K^N$ is reproducing,  we can rewrite this formula:
\begin{align} \notag
\Cu^n_{K^N}[f]
&= - \sum_{\k \vdash n}   \frac{(-1)^{\ell(\k)}}{\ell(\k)} \M(\k) 
\underset{x_{0} =x_n}{\int_{\X^n}} f(x_1)^{k_1} \cdots  f(x_l)^{k_l} \prod_{1\le j\le n} K^N(x_j, x_{j-1}) d\mu^n(x) \\
&\label{cumulant_0}
= -\underset{x_{0} =x_n}{\int_{\X^n}}\Upsilon^n_0[f](x)  \prod_{1\le j\le n} K^N(x_j, x_{j-1}) d\mu^n(x)  ,
\end{align}
where for any $\k\in\mho$,
\begin{equation}\label{Upsilon_0}
\Upsilon_0(\k) = \frac{(-1)^{\ell(\k)}}{\ell(\k)} \M(\k)  . 
\end{equation}
A simple observation which turns out to be very important when it comes to asymptotics is  that, by lemma~\ref{thm:combinatorics}, for all $n\ge 2$,
\begin{equation} \label{normalization_0}
\sum_{\k \vdash n} \Upsilon_0(\k)  =\gamma_0^n =0 . 
\end{equation}

For any $m\in\N$, we define\footnote{By convention $ {l \choose m} = 0 $ if $m>l$.}
\begin{equation} \label{Upsilon}
  \begin{cases}
\displaystyle  \Upsilon_m(\k) = \frac{(-1)^{\ell(\k)}}{\ell(\k)} {\ell(\k) \choose m} \M(\k)  & \text{if } \ell(\k) \ge 2 \\
\displaystyle  \Upsilon_m(\k) =- \delta_{1}(m)-\gamma^n_m  & \text{if } \k=n
  \end{cases} .
  \end{equation}
These functions are constructed so that we have for all $n,m\in\N$, 
\begin{equation} \label{normalization}
 \sum_{\k \vdash n} \Upsilon_m(\k) = 0 .
 \end{equation}

According to formula \eqref{Soshnikov}, the cumulants of the process $\widehat{\Xi}$ with correlation kernel $\widehat{K}_N = p_N K^N$ are given by 
  \begin{equation*}
\Cu^n_{\widehat{K}^N}[f]
= % \sum_{\k \vdash n}   \frac{(-1)^{\ell(\k)}}{\ell(\k)} \M(\k) 
- \sum_{l=1}^n   \frac{(-1)^l}{l} p_N^l \sum_{\begin{subarray}{c} \k \vdash n \\ \ell(\k) =l \end{subarray}} \M(\k)
\underset{x_{0} =x_l}{\int_{\X^l}} f(x_1)^{k_1} \cdots  f(x_l)^{k_l} \prod_{1\le j\le l} K^N(x_j, x_{j-1}) d\mu^l(x) . 
\end{equation*}
Since the kernel $K^N$ is reproducing, if we set $q_N=1-p_N$, using the binomial formula, we obtain
\begin{equation} \label{cumulant}
    \begin{aligned}  
\Cu^n_{\widehat{K}^N}[f]
= \Cu^n_{K^N}[f]  %& + q  \bigg\{ \gamma^n_1 \int_\X  f(x)^{n} K(x,x) d\mu(x)  +  \underset{x_{0} =x_n}{\int_{\X^n}}     \Upsilon^n_1[f](x)  \prod_{j\le n} K(x_j, x_{j-1}) d\mu^n(x)  \bigg\}\\
&- \sum_{m=1}^n (-q_N)^m  \gamma^n_m \int_\X  f(x)^{n} K^N(x,x) d\mu(x)  \\
&-    \sum_{m=1}^n (-q_N)^m  \underset{x_{0} =x_n}{\int_{\X^n}}     \Upsilon^n_m[f](x)  \prod_{1\le j\le n} K^N(x_j, x_{j-1}) d\mu^n(x) .
\end{aligned}
\end{equation}

 We are now ready to state our general result from which Theorems~\ref{thm:crossover_2},  \ref{thm:crossover_3}, \ref{thm:crossover_1} and \ref{thm:crossover_4} in the introduction follow. 
 
 \begin{theorem} \label{thm:main}
 Let $0<q_N<1$ be a sequence which converges to 0 as $N\to+\infty$.
 Under our general assumptions above, let $f_N$ be a sequence of functions for which
%and $L_N>0$ be a sequence which is either $1$ or goes to $+\infty$ as $N\to+\infty$. 
 %Let $f:\X \to \R$ be a locally integrable function and $f_N = f\big(L_N(\cdot-x_0)\big)$ for a fixed $x_0 \in \X$. 
 the cumulants $\Cu^n_{K^N}[f_N]$ are  well-defined for all $n, N \in \N$ and the following conditions hold:
  \begin{enumerate}
  \item There exists a (Radon) measure $\eta$ on $\X$, a function $f\in L^p(\eta)$ for any $p\ge 2$,  and a   sequence $M_N \nearrow +\infty$ as $N\to+\infty$ such that for all $n\ge 1$, 
    \begin{equation} \label{ass1}
\frac{1}{M_N} \int_\X f_N(x)^n  K^N(x,x) d\mu(x) \simeq   \int_\X f(x)^n d\eta(x). 
  \end{equation}
\item  For all $n\ge 2$ and all  $m\ge 1$, as $N\to+\infty$
\begin{equation}\label{ass2}
 \left| \underset{x_{0} =x_n}{\int_{\X^n}}     \Upsilon^n_m[f_N](x)  \prod_{1\le j\le n} K^N(x_j, x_{j-1}) d\mu^n(x) \right| = o\left(q_N^{-1} \vee M_N \right). 
\end{equation}
\item There exists $\sigma >0$ such that for all $n\in \N$, 
\begin{equation} \label{CLT0}
\lim_{N\to\infty} \Cu^n_{K^N}[f_N] = \begin{cases} \sigma^2 &\text{if } n=2 \\ 0 &\text{if } n>2 
\end{cases} .
\end{equation}
\end{enumerate}
Then, depending on the parameter  $T_N =  q_N M_N >0$,  we distinguish three different asymptotic regimes for the linear statistic $\widehat{\Xi}(f_N)$ of the   thinned point process with density $p_N=1-q_N$:
  \begin{enumerate}
  \item[i)] If $T_N \to 0$ as $N\to+\infty$, 
 \[  \widehat{\Xi}(f_N) - \E\big[\widehat{\Xi}(f_N) \big]  \Rightarrow \No(0,\sigma^2) . \]
\item[ii)] If $T_N \to \tau$ with $\tau>0$ as  $N\to+\infty$,
\[  \widehat{\Xi}(f_N) - \E\big[\widehat{\Xi}(f_N) \big]  \Rightarrow \mathrm{X} + \Lambda_{\tau\eta}(f), \]
where $\mathrm{X} \sim \No(0,\sigma^2)$ and $\Lambda_{\tau\eta}$ is Poisson process on $\X$ with intensity $\tau\eta$  independent from $\mathrm{X}$. 
\item[iii)] If $T_N \to +\infty$ as $N\to+\infty$, 
 \[  \frac{\widehat{\Xi}(f_N) - \E\big[\widehat{\Xi}(f_N) \big]}{\sqrt{T_N}}  \Rightarrow \No\Big(0, \int f(x)^2 d\eta(x) \Big) . \]
 \end{enumerate}
 \end{theorem}

 Before giving our proof of theorem~\ref{thm:main}, let us make two remarks about the assumptions \eqref{ass1}--\eqref{CLT0}.

 \begin{remark}
For any function $g:\X\to\R_+$,  we have  $\displaystyle\E[\Xi(g)] = \int_\X g(x)  K^N(x,x) d\mu(x)$, so that one can interpret \eqref{ass1} as a condition about the mean of the point process $\Xi$. 
In contrast, \eqref{ass2} can be seen as a  condition about the fluctuations of the incomplete point process.
We also implicitely assume that for $n=2$, the RHS of \eqref{ass1} is positive, so that the measure $\eta$ is non-trivial and  puts mass on the support of $f$. 
 Then the random variables $\No\big(0, \int f(x)^2 d\eta(x) \big)$ and $\Lambda_{\tau\eta}(f)$ are non zero for any $\tau>0$.  In order to handle mesoscopic linear statistics, we have allowed our test functions to depend on the parameter $N$. However, for simplicity, one can think of the case where $f_N = f$ is a smooth and compactly supported test function.
%Finally, let us emphasize that in the conditions \eqref{ass1} and \eqref{ass2}, we do not require any uniformity in the parameters $n,m \ge 1$.   
\end{remark}

 \begin{remark}
Instead of \eqref{CLT0}, we could assume that the cumulants of the linear statistics $\Xi(f_N)$ converge to that of a random variable $\mathrm{X}$ which is not necessarily Gaussian. Then, the conclusion $ii)$ of theorem~\ref{thm:main} remains true and we obtain a crossover from a non-Gaussian process to a Poisson process. For instance, this more general situation arises when considering linear statistics of 1-dimensional log-gases in the multi-cut regime.    
 \end{remark}

 \begin{proof}
 Observe that it follows from formula \eqref{cumulant} and the condition \eqref{ass2} that 
 the cumulants of the linear statistic $\widehat{\Xi}(f_N)$ satisfy for all $n\ge 2$
  as $N\to+\infty$, 
   \[ 
\Cu^n_{\widehat{K}^N}[f_N]
= \Cu^n_{K^N}[f_N]
- \sum_{m=1}^n (-q_N)^m  \gamma^n_m \int_\X  f_N(x)^{n} K^N(x,x) d\mu(x)  + o(1 \vee T_N) ,
\]
where $T_N =  q_N L_N$.  
Then, using the condition  \eqref{ass1}, we obtain for any  $n\ge 2$, as $N\to+\infty$, 
 \begin{equation} \label{limit4}
\Cu^n_{\widehat{K}^N}[f_N]
= \Cu^n_{K^N}[f_N] + T_N  \sum_{m=0}^{n-1} (-q_N)^m  \gamma^n_{m+1} \left(  \int_\X  f_N(x)^{n}  d\eta(x)  + o(1) \right)  + o(1 \vee T_N) .
\end{equation}
Let us observe that in the previous sum, regardless of the regime we consider, since $q_N\to0$ as $N\to+\infty$, only the term $m=0$ is asymptotically relevant.
For instance, if we assume that  $T_N = \tau +o(1)$ with $\tau \ge 0$, this implies that
   \[ 
\Cu^n_{\widehat{K}^N}[f_N]
= \Cu^n_{K^N}[f_N] + \tau  \gamma^n_{1}  \int_\X  f_N(x)^{n}  d\eta(x)+ o(1)  . 
\]
On the one hand by Lemma~\ref{thm:combinatorics}, since $\gamma^n_1 = (-1)^n$ and using the condition \eqref{CLT0}, we obtain for any $n\ge 2$
 \begin{equation}  \label{limit3} 
\lim_{N\to+\infty} \Cu^n_{\widehat{K}^N}[f_N]
=  \sigma^2 \1_{n = 2} + (-1)^n \tau   \int_\X  f_N(x)^{n}  d\eta(x) . 
\end{equation}
In the regime $i)$ -- which corresponds to $\tau=0$ -- this shows that the linear statistic $\widehat{\Xi}(f_N)$, once centered, converges in distribution (as well as in the sense of moments) to a Gaussian random variable with variance $\sigma^2$. 
 Let us observe that by formula \eqref{Poisson}, if $\tau>0$, the second term on the RHS of \eqref{limit3}
 corresponds to the $n^{\rm th}$ cumulant of the random variable $\Lambda_{\tau \eta}(f)$. This proves the claim in the regime $ii)$. 
 
 \medskip
 
 On the other hand, in the regime $iii)$ where $T_N\to+\infty$, we see from \eqref{limit4} that the variance 
$\Cu^n_{\widehat{K}^N}[f_N]$ diverges as $N\to+\infty$. Thus, in order to have a non-trivial limit, we need to renormalize the linear statistic $\widehat{\Xi}(f_N)$. Namely, we consider instead the test function $g_N = f_N /\sqrt{T_N}$ and it follows from \eqref{limit4} that  for all $n\ge 2$, as $N\to+\infty$
 \[
 \Cu^n_{\widehat{K}^N}[g_N] =   \1_{n=2}\left(  \int_\X  f(x)^{2}  d\eta(x)  + o(q_N) \right)  + o(1). 
 \]
 These asymptotics show that in the  the regime $iii)$, $\frac{\widehat{\Xi}(f_N) - \E\big[\widehat{\Xi}(f_N) \big]}{\sqrt{T_N}}$ converges in distribution (as well as in the sense of moments) to a centered Gaussian random variable with variance $\int_\X  f(x)^{2}  d\eta(x)$. 
 \end{proof}

\section{Transition for Coulomb gases in two dimensions}\label{sect:2}

\subsection{Asymptotics of the correlation kernel} \label{sect:asymp}

In this section, we begin by reviewing the basics of the theory of eigenvalues of random normal matrices developed  by Ameur, Hedenmalm and Makarov. In particular, we are interested in the properties of the correlation kernel \eqref{kernel} in the bulk of the gas.  
% We let $$ \S_V := \big\{z\in\C : \Delta V (z) >0 \big\} \cap \overline{\S_V}_V^\circ .$$
%
%
It has been established in \cite{AHM11} that, if the potential $V$ is real-analytic, the equilibrium measure  is $\varrho_V  = 2\Delta V \1_{\S_V}$ and the droplet $\overline{\S_V}$ is a compact set with a  nice boundary. 
Moreover, in order to compute the asymptotics of the cumulants of a smooth linear statistic, instead of working with  the correlation kernel $K^N_V$, one can use the so-called approximate Bergman kernel:
 \begin{equation}\label{B_kernel}
B^N_V(z,w) = \big( N b_0(z, \bar{w}) + b_1(z, \bar{w}) \big) e^{N \{ 2 \Phi(z, \bar{w}) - V(z) - V(w) \} } . 
\end{equation}
The functions  $b_0(z,w)$, $b_1(z,w)$ and $\Phi(z,w)$ are the (unique) bi-holomorphic functions defined in a neighborhood  in $\C^2$ of the set $\big\{ (z ,\bar{z}) : z\in \S_V \big\}$  such that 
$b_0(z,\bar{z}) = 2 \Delta V(z)$, $b_1(z,\bar{z}) = \frac{1}{2} \Delta \log( \Delta V)(z)$, and $\Phi(z,\bar{z}) = V(z)$.

\begin{lemma}[Lemma~1.2 in \cite{AHM11}, proved in  \cite{Berman09, AHM10}] \label{thm:Bergman_kernel}
For any $x_0 \in \S_V $, there exists $\epsilon_0>0$ and $C_0>0$ so that when the dimension $N$ is sufficiently large, we have for all $z, w  \in \D(x_0 , \epsilon_0)$, 
\begin{equation} \label{approximation_1}
\left| K^N_V(z, w) - B^N_V(z,w) \right| \le C_0 N^{-1} . 
\end{equation}
\end{lemma} 
 
Moreover, at sufficiently small mesoscopic scale, up to a gauge transform, the asymptotics of the  approximate Bergman kernel $B^N_V$ is universal.

\begin{lemma} \label{thm:Ginibre_approximation} 
Let $\kappa>0$ and  $\epsilon_N= \kappa N^{-1/2}\log N $ for all $N\in \N$.
For any  $x_0 \in \S_V$, there exists $\varepsilon_0>0$ and  a function $\h : \D(0,\varepsilon_0)  \to \R$ such that if  the parameter $N$ is sufficiently large, the function 
\begin{equation} \label{B}
\widetilde{B}^N_{V, x_0}( u ,v)  =\frac{B^N_V(x_0 + u , x_0 + v)e^{i N \h(u)}}{e^{i N \h(v)}} 
\end{equation} 
satisfies
\begin{equation} \label{asymptotics_B}
\widetilde{B}^N_{V, x_0}( u , v) = K^\infty_{N \varrho_V(x_0)}(u, v)\left\{ 1 + \underset{N\to\infty}{O}\big(  (\log N)^{2}\epsilon_N   \big)  \right\}
\end{equation}
uniformly for all $u, v \in \D(0, \epsilon_N)$, where $K^\infty_{N \varrho_V(x_0)}$ is the \Ginibre kernel with density $N \varrho_V(x_0)= 2N \Delta V(x_0)$. 
\end{lemma}

A key ingredient in the paper \cite{AHM11}, as well as \cite{RV07b},  is to reduce the domain of integration in formula \eqref{cumulant_0}, using the exponential  off-diagonal decay of the correlation kernels $K_V^N$, to a set where we can use the asymptotics \eqref{asymptotics_B} -- see the next lemma.
For completeness, the proofs of lemmas~\ref{thm:Ginibre_approximation} and~\ref{thm:localization} are given in the appendix~\ref{A:damping}. 

\begin{lemma} \label{thm:localization}
Let $n\in\N$ and $\epsilon_N = \kappa N^{-1/2} \log N$ for some constant $\kappa>0$ which is sufficiently large compared to $n$. We let
 \begin{equation} \label{A}
\mathscr{A}(z_0; \epsilon)
= \big\{ \z\in\C^n  :  |z_j - z_{j-1}| \le \epsilon \text{ for all } j=1, \dots n  \big\} .  
\end{equation}
Let $\S$ be a compact subset of $\S_V$, $N_0\in\N$,  and $F_N : \C^{n+1} \to \R $ be a sequence of  continuous functions such that 
\begin{equation}\label{estimate_0}
 \sup\big\{ |F_N(z_0, \z)| : \z \in \C^{n} , N  \ge N_0 \big\}  \le C \1_{z_0 \in\S} .
 \end{equation}
We have
\begin{align}
&\label{localization} \underset{z_0 = z_{n+1}}{\int_{\C^{n+1}}}\hspace{-.2cm} F_N(z_0, \z) \prod_{j=0}^n K^N_V(z_{j}, z_{j+1}) d\A(z_0) d\A^{n}(z)  \\
&\notag\hspace{1cm}
= \int_\S d\A(z_{0})  \underset{z_{n+1} = z_{0}}{\int_{\mathscr{A}(z_{0}; \epsilon_N)}}\hspace{-.2cm} F_N(z_0,\z)  \prod_{j=0}^n K^N_V(z_{j}, z_{j+1}) d\A^{n}(\z)\ +  \O(N^{-1}) .
\end{align}
\end{lemma}

\begin{remark}  \label{rk:localization}
Recall  that $K^\infty_\rho$, \eqref{Ginibre}, denotes the correlation kernel of the \Ginibre process with density $\rho>0$. 
Using the fact that for all $z, w\in\C$, 
\begin{equation} \label{Ginibre_bound}
\big| K^\infty_\rho(z,w)\big| = \rho e^{-\rho |z-w|^2/2} , 
\end{equation}
it is easy to obtain the counterpart of Lemma~\ref{thm:localization} for the \Ginibre kernel, see lemma~\ref{thm:delocalization}. 
\end{remark}

\subsection{Proof of theorem~\ref{thm:crossover_2}} \label{sect:global2}

In this section, we show how to apply theorem~\ref{thm:main} for Coulomb gases in the global regime by relying on the asymptotics from section~\ref{sect:asymp} for the correlation kernel $K^N_V$. 
First, observe that  with $f_N = f$ and $M_N = N$,  the  assumptions \eqref{ass1}  follow immediately from the law of large numbers \eqref{LLN}. 
Then the measure $d\eta = \varrho_V d\mu$ is absolutely with respect to $\mu$ with compact support. 
Moreover, if $f\in \mathcal{C}_c^3(\S_V)$, then the conditions \eqref{CLT0} are well-known from \cite[Theorem~4.4]{AHM11} with $\sigma^2 = \E[\mathrm{X}(f)]$ according to formula \eqref{H1_noise}. So, our main technical challenge  is to obtain the estimates \eqref{ass2} for a large class of test functions.
The first step of the proof is the following approximation.

%in order  to complete the proof of theorem~\ref{thm:crossover_2}, we only need to show that the 2-dimensional log-gases satisfy the condition \eqref{ass2}. 
%The technical step of the proof is the next proposition. In fact, if we choose $L_N=1$ and $x_0=0$, the proof of proposition~\ref{thm:approximation} still works as long as the support of the test function $\K \subset\S_V$.
%

%The  idea of the proof is to use Lemma~\ref{thm:localization} and then to  approximate the function $\Upsilon^n_0[f]$ by a multivariate polynomial of degree 2 using Taylor's theorem.

\begin{proposition} \label{thm:approximation}
Let $\K\subset \C$ be a compact set, $f\in \mathcal{C}^3_c(\K)$ and let $\Upsilon: \mho \to\R$ be any map such that  $\sum_{\k \vdash n} \Upsilon(\k) = 0$ for all $n\ge2$. Let  $L_N $ be an increasing  sequence such that $L_N^{-1} (\log N)^4 = o(1)$ and $L_N  = o \big(\sqrt{N}/(\log N)^{3}\big)$ as $N\to+\infty$. We also denote by $H^n(\lambda; \w )$  the second order Taylor polynomial at $0$ of the function 
$\w\in\C^n \mapsto   \Upsilon^{n+1}[f](\lambda, \lambda+ \w)$. 
Fix  $x_0 \in \S_V$ and let $f_N(z)=f(L_N(z- x_0))$ as in \eqref{meso}.  
Then, we have for any  $n\ge 1$, as $N\to+\infty$, 
\begin{equation}
\begin{aligned} \label{approximation}
& \underset{z_0 = z_{n+1}}{\int_{\C^{n+1}}}\hspace{-.2cm}  \Upsilon^{n+1}[f_N](z_0, \z) \prod_{j=0}^n K^N_V(z_{j}, z_{j+1}) d\A(z_0) d\A^{n}(\z)  \\
&\qquad\qquad \simeq \int_{\mathcal{K}} d\A(\lambda)
\hspace{-.3cm}\underset{w_0=w_{n+1}=0}{\int_{\C^n}}\hspace{-.4cm}
H^n(\lambda; \w)
\prod_{j=0}^{n} K^\infty_{\eta_N(\lambda)}(w_{j}, w_{j+1}) d\A^n(\w)  , 
\end{aligned} 
\end{equation}
where the density is given by $\eta_N(\lambda) = N L_N^{-2} \varrho_V(x_0+ \lambda/L_N)$. %In particular, note that by assumptions $\eta_N(\lambda) \to\infty$ for all $\lambda\in\S_0$. 
\end{proposition}

\begin{remark}\label{rk:scales}
Observe that under the assumptions of proposition~\ref{thm:approximation}, we have $\eta_N(\lambda) \to +\infty$ as $N\to+\infty$ uniformly for all $\lambda \in \K$. 
Moreover, it follows from the proof below that in the global regime where $L_N =1$ and  $x_0 = 0$, provided that  $\K \subset \S_V$, the estimates \eqref{approximation} remain valid with an extra error. 
Namely, we obtain  for all $n\ge 1$,
\begin{equation}\begin{aligned} \label{approximation'}
& \underset{z_0 = z_{n+1}}{\int_{\C^{n+1}}}\hspace{-.2cm}  \Upsilon^{n+1}[f_N](z_0, \z) \prod_{j=0}^n K^N_V(z_{j}, z_{j+1}) d\A(z_0) d\A^{n}(\z)  \\
&\qquad
= \int_{\mathcal{K}} d\A(\lambda)
\hspace{-.3cm}\underset{w_0=w_{n+1}=0}{\int_{\C^n}}\hspace{-.4cm}
H^n(\lambda; \w)
\prod_{j=0}^{n} K^\infty_{ N \varrho_V(\lambda)}(w_{j}, w_{j+1}) \  d\A^n(\w)  
+\O\left((\log N)^4 \right) .
\end{aligned} 
\end{equation}
\end{remark}

\proof We let $F_N = \Upsilon^{n+1}[f_N]$ and 
$$\mathscr{J}^{n}_N : =  \underset{z_{n+1}=z_0}{\int_{\C^{n+1}}}\hspace{-.2cm} F_N(z_0, \z) \prod_{j=0}^n K^N_V(z_{j}, z_{j+1}) d\A(z_0) d\A^{n}(\z)  . $$
%Note that we abuse of notation in the previous formula and denoted $F_N( \z) = F_N(z_0, z_1,\dots, z_n)$. 
Since $x_0\in\S_V$, there exists a compact set $\S\subset \S_V$ so that $\supp(f_N) \subseteq \S$ when the parameter $N$ is sufficiently large. Then, according to formula \eqref{Upsilon_f}, the function $F_N$ satisfies the assumption \eqref{estimate_0}.  Thus, by lemma~\ref{thm:localization}, we obtain as $N\to+\infty$
\begin{equation} \label{cumulant_1}
\mathscr{J}^{n}_N \simeq \int_{\S} d\A(z_0)  \underset{z_{n+1} = z_0}{\int_{\mathscr{A}(z_0; \epsilon_N)}}\hspace{-.3cm} F_N(z_0,\z)  \prod_{j=0}^n K^N_V(z_{j}, z_{j+1}) d\A^{n}(\z) .
\end{equation}  
By \eqref{A}, the set
\begin{equation*}
\mathscr{A}(z_0; \epsilon)
\subset \big\{ \z\in\C^n :   z_1, \dots, z_n \in  \D(z_0,n \epsilon_N) \big\} 
\end{equation*}
and we can apply lemma~\ref{thm:Bergman_kernel} to replace the kernels $K^N_V(z_j,z_{j+1})$ in formula \eqref{cumulant_1}.
 Namely, if $z_{n+1} = z_0 $ and $\z \in \mathscr{A}(z_0; \epsilon)$, then
 %Namely for any $z_0 \in \S_0$ and for all $\z \in \mathscr{A}(z_0;\epsilon_N)$,  if  $z_{n+1} = z_0$, then
\begin{equation*}
\Bigg| \prod_{j=0}^n K^N_V(z_{j}, z_{j+1}) - 
 \prod_{j=0}^{n} B^N_V(z_{j}, z_{j+1}) \Bigg| \le C \sum_{k=1}^{n+1} N^{-k} S_N^{n+1-k} ,  
\end{equation*}
where $S_N =  \sup\big\{ |B^N_V(z,w)|  :  z, w \in \D(z_0, n \epsilon_N)  , z_0 \in \S  \big\} $. 
By  lemma~\ref{thm:Ginibre_approximation}, we have for  all  $u,v \in \D(0, n\epsilon_N) $, 
$$
\big| B^N_V(z_0 + u , z_0 + v) \big| \le C
\big|  K^\infty_{N\varrho_V(z_0)}(u,v) \big| .
$$
and, by formula \eqref{Ginibre_bound}, this implies that $S_N  \le  C N$.
If we combine these estimates with formula \eqref{cumulant_1}, since the functions $F_N$ are uniformly bounded, we obtain
\begin{align*} 
\mathscr{J}^{n}_N
=  \int_{\S} d\A(z_0)  \underset{z_{n+1} = z_0}{\int_{\mathscr{A}(z_0; \epsilon_N)}}\hspace{-.3cm} F_N(z_0,\z)  \prod_{j=0}^n B^N_V(z_{j}, z_{j+1}) d\A^{n}(\z)
%\hspace{3cm}
+\O\left(  N^{n-1}  \int_{\S} d\A(z_0)  \big|\mathscr{A}(z_0; \epsilon_N)\big| \right) ,
\end{align*}  
where $|\mathscr{A}|$ denotes the Lebesgue measure of the set $\mathscr{A}$. 
By definition, $\epsilon_N = \kappa N^{-1/2} \log N$ so that 
$ \big| \mathscr{A}(z_0;\epsilon_N)  \big|  \le C N^{-n} (\log N)^{2n}$ 
for all $z_0\in \C$. Thus, the previous error term converges to 0 like $(\log N)^{2n}/ N$.
 Hence, if we make the change of variables $\z = z_0 +\u$ and the appropriate {\it gauge transform} in the previous integral, according to formula \eqref{B}, we obtain
\begin{equation}  \label{cumulant_2}
\mathscr{J}^{n}_N
\simeq  \int_{\S} d\A(z_0)  \underset{u_{n+1} = u_0=0}{\int_{\mathscr{A}(0; \epsilon_N)}}\hspace{-.3cm} F_N(z_0 +\u) \prod_{j=0}^n \widetilde{B}^N_{V, z_0}(u_{j}, u_{j+1}) d\A^{n}(\u) . 
\end{equation}
Note that in formula \eqref{cumulant_2}, the integral  is over a small subset of the surface $\{ \u \in \C^{n+2} : u_0 = u_{n+1} =0\} $ and we denote 
$F_N(z_0 +\u) = F_N(z_0, z_0+u_1,\dots, z_0+u_n)$. 
%We will use this convention in the rest of the proof.
Then, we  can apply lemma~\ref{thm:Ginibre_approximation} to replace the kernel $ \widetilde{B}^N_{V, z_0}$ by $K^\infty_{N \varrho_V(z_0)}$ in formula \eqref{cumulant_2}, we obtain 
\begin{equation*}  
\mathscr{J}^{n}_N
\simeq  \int_{\S} d\A(z_0)  \underset{u_{n+1} = u_0=0}{\int_{\mathscr{A}(0; \epsilon_N)}}\hspace{-.3cm} F_N(z_0 +\u)\chi_N(z_0,\u)  \prod_{j=0}^n K^\infty_{N \varrho_V(z_0)}(u_{j}, u_{j+1})  d\A^{n}(\u)  ,
\end{equation*}
where $\displaystyle \chi_N(z_0,\u) = 1 + \underset{N\to\infty}{O}\big(  (\log N)^{2}\epsilon_N   \big)$ uniformly for all $z_0\in\S$ and all $\u\in \mathscr{A}(0; \epsilon_N)$. 

\medskip

Let $F=\Upsilon^{n+1}[f]$, $\delta_N = \epsilon_N L_N$ and $\eta_N(\lambda) =N L_N^{-2} \varrho_V(x_0 + \lambda/L_N) $. By definition, 
$F_N(z_0 +\u) = F\big(L_N(z_0- x_0 +\u )\big) $ and 
we can make the change of variables $\lambda = L_N(z_0-x_0)$ and $\w= L_N \u$ to get rid of the scale $L_N$ and $x_0$ in the previous integral. Using the obvious scaling property of the \Ginibre kernel, \eqref{Ginibre}, we obtain
\begin{equation} \label{cumulant_3}  
\mathscr{J}^{n}_N
\simeq  \int_{\K} d\A(\lambda)  \underset{w_{n+1} = w_0=0}{\int_{\mathscr{A}(0; \delta_N)}}\hspace{-.3cm} F(\lambda+\w)\widetilde{\chi}_N(\lambda,\w)  \prod_{j=0}^n K^\infty_{\eta_N(\lambda)}(w_{j}, w_{j+1})  d\A^{n}(\w) ,  
\end{equation}
where $\displaystyle \widetilde{\chi}_N(\lambda,\w)= 1 + \underset{N\to\infty}{O}\big(  (\log N)^{2}\epsilon_N   \big)$ uniformly for all $\lambda\in\K$ and for all $\u\in \mathscr{A}(0; \delta_N)$. 
Here we used that the test function $f$ is supported in the set $\K$. 
The condition $\sum_{\k \vdash n+1} \Upsilon(\k) = 0$ implies that $F(\lambda+0)=0$ for all $\lambda \in \C$ so that   for  all $\w\in \mathscr{A}(0; \delta_N)$, 
$$ \big| F(\lambda+\w) \big| \le C \delta_N .$$
Moreover, by formula~\eqref{Ginibre_bound}, we have  for any $n \in \N$, 
\begin{equation} \label{estimate_5}
\prod_{j=0}^{n}  \big|K^\infty_{\rho}(w_{j}, w_{j+1})\big| \le \rho^{n+1}   \prod_{j=1}^{n}  e^{- \rho | v_j |^2 /2 } 
\end{equation}
where $v_j = w_j - w_{j-1}$ for all $j=1,\dots, n$. Hence, we see that
$$
\bigg| \underset{w_{n+1} = w_0=0}{\int_{\mathscr{A}(0; \delta_N)}}\hspace{-.3cm} F(\lambda+\w) \prod_{j=0}^n K^\infty_{\eta_N(\lambda)}(w_{j}, w_{j+1})  d\A^{n}(\w)  \bigg| \le C \delta_N \eta_N(\lambda) 
$$
and, since $\eta_N(\lambda) \le C N L_N^{-2}$ for all $\lambda\in \K$, we deduce from formula \eqref{cumulant_3} that
\begin{align}  \label{cumulant_4}  
\mathscr{J}^{n}_N
= \int_{\K} d\A(\lambda)  \underset{w_{n+1} = w_0=0}{\int_{\mathscr{A}(0; \delta_N)}}\hspace{-.3cm} F(\lambda+\w) \prod_{j=0}^n K^\infty_{\eta_N(\lambda)}(w_{j}, w_{j+1})  d\A^{n}(\w)  \\
+ \notag\O\big( N L_N^{-2} \delta_N \epsilon_N (\log N)^2 \big) .
\end{align}
Recall that $\delta_N= L_N \epsilon_N$ and $\epsilon_N=\kappa N^{-1/2}\log N$, so that the error term in \eqref{cumulant_4} is of order $(\log N)^4 L_N^{-1}$. 
Moreover,  if  $L_N = o\big( \sqrt{N}/\log N \big)$,   a Taylor approximation shows that for any $\w \in \mathscr{A}(0; \delta_N) $, 
\begin{equation*} \label{Taylor}
F(\lambda,\lambda+w_1,\dots,  \lambda+w_n)  =H^n( \lambda; \w) + \O(\delta_N^3) . 
\end{equation*}
Using the estimate \eqref{estimate_5} once more, by formula \eqref{cumulant_4}, this implies that
\begin{align}  \label{cumulant_5}  
\mathscr{J}^{n}_N
=  \int_{\K} d\A(\lambda)  \underset{w_{n+1} = w_0=0}{\int_{\mathscr{A}(0; \delta_N)}}\hspace{-.3cm} H^n( \lambda; \w) \prod_{j=0}^n K^\infty_{\eta_N(\lambda)}(w_{j}, w_{j+1})  d\A^{n}(\w)  \\
+ \notag\O\big( N L_N^{-2} \delta_N^3 \vee (\log N)^4 L_N^{-1} \big) .
\end{align}
By lemma~\ref{thm:delocalization}, the leading term in formula \eqref{cumulant_5} has the same limit (up to an arbitrary small error term) as
$$
\int_{\K} d\A(\lambda)  \underset{w_{n+1} = w_0=0}{\int_{\C^n}}\hspace{-.3cm} H^n( \lambda; \w) \prod_{j=0}^n K^\infty_{\eta_N(\lambda)}(w_{j}, w_{j+1})  d\A^{n}(\w) 
$$
and, since  $N L_N^{-2} \delta_N^3 \to0$ when $L_N = o \big( \sqrt{N}/(\log N)^3 \big)$, this completes the proof.  \qed\\

%Let $0<q_N<1$ be a sequence which converges to 0 as $N\to\infty$. Let $p_N = 1- q_N$ and denote $\widehat{K}^N = p_NK^N$.   

%By formula \eqref{cumulant}, proposition~\ref{thm:approximation} shows that the asymptotics of the cumulants of linear statistics of the  log-gases in dimension 2 is reduced to compute the integral of a multivariate polynomial of degree 2 against a product of Ginibre kernels.  
%This analysis has been performed in \cite{RV} and since our setting is a bit different, we will review the main steps in the rest of this section. \\

Since the function $\w\mapsto H^n(\lambda; \w )$ is a multivariate polynomial of degree 2,  the leading term in the asymptotics~\eqref{approximation} can be computed explicitly using the reproducing property of the \Ginibre kernel; see for instance~\cite{RV07b}. 
For any $\rho>0$, the function  $(z,w) \mapsto e^{\rho z \bar{w}}$ is the reproducing kernel for the Bergman space with weight  $\rho e^{-\rho|z|^2/2}$ on $\C$. This implies that for any $w_1, w_2 \in \C$ and for all integer $k\ge 0$,
\begin{equation*} %\label{reproducing}
\begin{cases}
\displaystyle \int_\C K^\infty_{\rho}(w_1,w_2) w_2^k K^\infty_{\rho}(w_2,w_3) d\A(w_2)  = w_1^kK^\infty_{\rho}(w_1,w_3)  \vspace{.2cm}   \\
%$$
%and taking complex conjugate:
%$$
\displaystyle \int_\C K^\infty_{\rho}(w_1,w_2) \bar{w_2}^k K^\infty_{\rho}(w_2,w_3) d\A(w_2)  = \bar{w_3}^kK^\infty_{\rho}(w_1,w_3)     
\end{cases} .
\end{equation*}

As a basic application  of these identities, we obtain the following lemma.

\begin{lemma} \label{thm:reproducing}
Let $n \ge 1$ and  $\rho>0$. For any polynomial $H(\w)$ of degree at most~2 in the variables $w_1,\dots, w_n, \bar{w_1}, \dots, \bar{w_n}$ such that $H(0)=0$, we have
\begin{equation} \label{reproducing}
 \underset{w_{n+1} = w_0=0}{\int_{\C^n}}\hspace{-.3cm} H^n( \w) \prod_{j=0}^nK^\infty_{\rho}(w_{j}, w_{j+1}) \  d\A^n(\w)  
%& =  \int_{\C^n} \mho[H](w) \prod_{j=0}^nK^\infty_{\rho}(w_{j}, w_{j+1}) \  d\A^n(w)  \\
  =    \sum_{1 \le r \le s\le n} \p_s \pb_r H |_{\w=0} . 
 \end{equation}
 %where 
 %$$ \mho[H](w) = \sum_{r, s =1}^n \p_r \pb_s H(0) w_r  \bar{w_s}  $$
%is the polynomial constructed from $H$ by keeping only the mixed terms.
\end{lemma}

Under the assumptions of proposition~\ref{thm:approximation}, since $\eta_N(\lambda) \to +\infty$ as $N\to+\infty$ uniformly for all $\lambda \in \K$, we deduce from lemma~\ref{thm:reproducing} that for any test function $f\in \mathcal{C}^3_c(\K)$, we have for all integers $n\ge 1$ and $m\ge 0$,  as $N\to+\infty$, 
 \begin{equation} \label{cumulant_6}
\underset{z_0 = z_{n+1}}{\int_{\C^{n+1}}}\hspace{-.2cm}  \Upsilon^{n+1}_m[f_N](z_0, \z) \prod_{j=0}^n K^N_V(z_{j}, z_{j+1}) d\A(z_0) d\A^{n}(\z)  
\simeq  \hspace{-.3cm} \sum_{2 \le r \le s \le n+1}\int_{\mathcal{K}}  \p_s \pb_r \Upsilon^{n+1}_m[f] (\lambda,\dots, \lambda)  d\A(\lambda) .
\end{equation}
Here we used that according to formulae \eqref{normalization_0} and \eqref{normalization}, we have for any $m \ge 0$ and $n\ge 1$, 
\[
\Upsilon^{n+1}_m[f](\lambda, \dots , \lambda)= f(\lambda)^{n+1} \sum_{\k \vdash n+1} \Upsilon_m(\k) = 0 .
\] 

In the macroscopic regime  ($L_N=1$, $x_0=0$ and $\K = \supp(f) \subset \S_V$), by remark~\ref{rk:scales}, this also shows that for any $n,m \ge 1$, 
\[
 \underset{z_0 = z_{n+1}}{\int_{\C^{n+1}}}\hspace{-.2cm}  \Upsilon^{n+1}_m[f](\z) \prod_{j=0}^n K^N_V(z_{j}, z_{j+1}) d\A^{n+1}(\z) = \O\left((\log N)^4 \right) .
\]

This shows that the estimate \eqref{ass2} with $M_N = N$ holds for any sequence $q_N \searrow 0$ as $N\to+\infty$. By theorem~\ref{thm:main}, this completes the proof of theorem~\ref{thm:crossover_2}.

\medskip

\subsection{Mesoscopic fluctuations for 2-dimensional Coulomb gases and the proofs of theorem~\ref{meso:2} and theorem~\ref{thm:crossover_3}}

In the mesoscopic regime, we claim  that the asymptotics \eqref{cumulant_6} with $m=0$ implies the central limit theorem~\ref{meso:2}.  Indeed, the fact that the cumulants of order $n\ge 3$ vanish in the large $N$ limit comes from the following combinatorial lemma.

\begin{lemma}[\cite{RV07b}, Lemma~9] \label{thm:RV}
 For any $n\ge 1$, let
\begin{align*}
\mathscr{Y}_n=
-\sum_{k\vdash n}
\Upsilon_0(\k) 
\left\{  \sum_{2 \le r < s\le \ell(\k)}  k_r k_s +   \sum_{r=2}^{\ell(\k)} k_r (k_r-n) \right\}.
  \end{align*}
  We have $\displaystyle
\mathscr{Y}_n = \begin{cases} 1 &\text{if } n=2 \\ 0 &\text{else} \end{cases} .
$
\end{lemma}

\begin{proof}[Proof of theorem~\ref{meso:2}] Let $\lambda\in\C$ and $\boldsymbol\lambda =(\lambda,\dots, \lambda)\in\C^{n+1}$. 
According to formula \eqref{Upsilon_f}, an elementary computation shows that for any $2 \le r < s\le n+1$,
\begin{equation} \label{combinatorics_1}
 \p_s \pb_r \Upsilon_0^{n+1}[f](\boldsymbol\lambda)
=   \p f(\lambda) \pb f(\lambda) f(\lambda)^{n-1} 
\sum_{\k \vdash n+1} \Upsilon_0(\k)  k_r k_s \1_{s \le \ell(k)}
 \end{equation}
 and
 \begin{align} \notag
 \p_r \pb_r \Upsilon_0^{n+1}[f]((\boldsymbol\lambda)
&=    \p f(\lambda) \pb f(\lambda) f(\lambda)^{n-1} 
\sum_{\k \vdash n+1} \Upsilon_0(\k)  k_r(k_r-1) \1_{r \le \ell(k)}\\
&\hspace{.2cm} \label{combinatorics_2}
+  \Delta f(\lambda)  f(\lambda)^{n}\sum_{\k \vdash n+1} \Upsilon_0(\k)  k_r \1_{r \le \ell(k)} .
 \end{align}
Since, by integration by parts, 
 $$
 \int_\C  \Delta f(\lambda)  f(\lambda)^{n} d\A(\lambda)
 = - n  \int_\C  \p f(\lambda) \pb f(\lambda) f(\lambda)^{n-1} d\A(\lambda) ,
 $$
we deduce from formulae \eqref{combinatorics_1} and \eqref{combinatorics_2} that
\begin{align*}
\sum_{2 \le r \le s \le n+1}\int_{\C}  \p_s \pb_r \Upsilon^{n+1}_m[f] (\boldsymbol\lambda)  d\A(\lambda) = \mathscr{Y}_{n+1}
\int_{\C} \p f(\lambda) \pb f(\lambda) f(\lambda)^{n-1} d\A(\lambda) . 
  \end{align*}
When $ L_N =N^\alpha$ and  $0<\alpha<1/2$, formulae \eqref{cumulant_0} and \eqref{cumulant_6} with $m=0$  imply  that for any $n\ge 1$,
\begin{equation}\label{cumulant_7}
\lim_{N\to\infty}\Cu^{n+1}_{K^N_V}[f_N] = \mathscr{Y}_{n+1}
\int_{\C} \p f(\lambda) \pb f(\lambda) f(\lambda)^{n-1} d\A(\lambda) . 
\end{equation}
By lemma~\ref{thm:RV}, this proves that for any test function $f\in\mathcal{C}^3_0(\C)$ and any $n\ge 2$, 
\begin{equation} \label{cumulant_9}
\lim_{N\to\infty}\Cu^{n}_{K^N_V}[f_N] = \begin{cases}
 \|f \|_{H^1(\C)}^2 &\text{if }n=2  \\
 0 &\text{else}
\end{cases}. 
\end{equation}
This shows that the centered mesoscopic linear statistics 
$\Xi(f_N) -\E\big[\Xi(f_N) \big] $ converges in distribution  as $N\to\infty$ to the mean-zero Gaussian random variable $\mathrm{X}(f)$.
\end{proof}

We are now ready to finish the proof of theorem~\ref{thm:crossover_3}. By lemma~\ref{thm:Bergman_kernel},  for any bounded function $f$ with compact support, we have for any $n\ge 1$, 
\begin{align*} 
 \int_\C f_N(z)^n  K^N_V(z,z) d\A(z)
 &= N  \int_\C f_N(z)^n 2 \Delta V(z) d\A(z) + \O(1) \\
 & = N L_N^{-2} \varrho_V(x_0)  \int_\C f(z)^n d\A(z) + \O(N L_N^{-3} ) .
\end{align*}
Here we used that the potential $V$ is smooth and $\rho_V = 2 \Delta V >0$ on a small neighborhood of the point $x_0 \in \S_V$. 
This implies the assumption \eqref{ass1} with $M_N = N L_N^{1-2\alpha} \varrho_V(x_0)$ -- since the parameter $\alpha<1/2$,  $M_N \nearrow +\infty$ as $N\to+\infty$. 
As we already pointed out, the asymptotics \eqref{cumulant_6} yield the assumption  \eqref{ass2} with an error which is  $O(1)$. 
Finally, the assumption \eqref{CLT0} was proved just above -- see \eqref{cumulant_9}. So, by theorem~\ref{thm:main}, this completes the proof of theorem~\ref{thm:crossover_3}.

\section{Transition for 1-dimensional log-gases}\label{sect:1}

\subsection{Asymptotics of orthogonal polynomials} \label{sect:aop}

In this section, we begin by reviewing basic facts about the asymptotics of orthogonal polynomials which are required for the proofs of theorem~\ref{thm:crossover_1} and theorem~\ref{thm:crossover_4}.
A comprehensive reference for the results discussed in this section is the book of Deift \cite{Deift99}. We assume that the potential $V \in C^2(\R)$ is a function which  satisfies the condition~\eqref{potential} and  we let $\Xi$  and $\widehat{\Xi}$ be the determinantal processes with correlation kernels $K^N_V$ and $\widehat{K}^N_V = p_N K^N_V$ respectively.

\medskip

The proof of theorem~\ref{thm:crossover_1} relies on a combinatorial method introduced in \cite{BD17} which consists in using the  three-terms recurrence relation of the orthogonal polynomials $\{P_k \}_{k=0}^\infty$ with respect to the measure 
$d\mu_N = e^{- 2N V(x)}dx$ to compute the cumulants of polynomial linear statistics. 
For any $N\in \N$, there exists two sequences $a^N_k >0$ and $b^N_k \in \R$ such that the orthogonal polynomials $P_k$ in \eqref{OP} satisfy 
\begin{equation}\label{recurrence}
xP_k(x) = a^N_k P_{k+1}(x) + b^N_k P_k(x) +  a^N_{k-1} P_{k-1}(x)  .
\end{equation}

In particular, the completion $\mathscr{L}_N$ of the space of polynomials with respect to $L^2(\R,\mu_N)$ is isomorphic to $L^2(\N_0)$ and formula \eqref{recurrence} implies that the multiplication by $x$ on $\mathscr{L}_N$ is unitary equivalent to applying the  {\bf Jacobi matrix}
\begin{equation}\label{Jacobi}
\mathbf{J} := \begin{bmatrix}
 b_0^N & a_0^N & 0 & 0 & 0  \\
 a_0^N & b_1^N & a_1^N & 0 & 0  &\mathbf{0} \\
  0 &  a_1^N  &  b_2^N  &  a_2^N  & 0 \\ 
 0 & 0 &  a_2^N &  b_3^N &  a_3^N \\
 &\mathbf{0} &  & \ddots &\ddots&\ddots
\end{bmatrix} .
\end{equation}   
We also let $\boldsymbol\Pi_N$ be the orthogonal projection on $\operatorname{span}\{e_0, \dots, e_{N-1}\}$ acting on  $L^2(\N_0)$. 
The connection with eigenvalues statistics comes from the fact that for any polynomial
$Q$ and  for any composition $\k\vdash n$, one has
\begin{align*} \notag
&\underset{x_{0} =x_l}{\int_{\X^l}} Q(x_1)^{k_1} \cdots  Q(x_l)^{k_l} \prod_{0\le j\le l} K^N_V(x_j, x_{j-1}) dx_1 \cdots dx_l \\
& \hspace{2cm}= \tr\big[Q(\J)^{k_1}\boldsymbol\Pi_N\cdots Q(\J)^{k_l}\boldsymbol\Pi_N \big] \\
& \hspace{2cm}= \sum_{m=0}^{N-1} \sum_{\pi \in \Gamma_{m}^n }   \prod_{j=1}^l \1_{\pi(k_1+\cdots+ k_j) < N} \prod_{i=0}^{n-1}
Q(\J)_{\pi(i)\pi(i+1)} ,
\end{align*}
where $\mathcal{G}$ denotes the   adjacency graph of the matrix $Q(\J)$ and
$$
\Gamma_{m}^n = \big\{ \text{paths }\pi  \text{ on the graph }\mathcal{G}\text{ of length }n \text{ such that }\pi(0)=\pi(n)=m  \big\} .
$$
Given a path $\pi$ of length $n$ and a composition $\k\vdash n$, we let
\begin{equation}
 \Phi_\pi^N(\k) := \1_{\hspace{-.1cm}\displaystyle\max_{1\le j <\ell(\k)}\pi(k_1+\cdots+ k_j) \ge N } .
\end{equation}
Observe that
$$
\prod_{1\le j <\ell(\k)} \1_{\pi(k_1+\cdots+ k_j) < N}  = 1-  \Phi_\pi^N(\k) , 
$$
so that by formula \eqref{cumulant_0},  the cumulants of a polynomial linear statistics are given by 
\begin{equation} \label{cumulant_15}
\Cu^n_{K^N_V}[Q]
= - \sum_{m =0}^{N-1} \sum_{\pi \in \Gamma_{m}^n } \prod_{i=0}^{n-1}
Q(\J)_{\pi(i)\pi(i+1)}  \sum_{\k \vdash n}\Upsilon_0(\k)\big( 1-\Phi_\pi^N(\k)\big) . 
\end{equation}

By definitions, there exists a constant $M>0$ which only depends on the degree of $Q$ and $n$ so that  $\Phi_\pi^N=0 $ for any path $\pi \in \Gamma^n_m$ as long as $m < N - M$. 
Since $  \sum_{\k \vdash n}\Upsilon_0(\k) = 0$ for all $n \ge 2$, formula \eqref{cumulant_15} implies that  
\begin{equation*}
\Cu^n_{K^N_V}[Q]
=  \sum_{m =N -M}^{N-1} \sum_{\pi \in \Gamma_{m}^n } \prod_{i=0}^{n-1}
Q(\J)_{\pi(i)\pi(i+1)}  \sum_{\k \vdash n}\Upsilon_0(\k)\Phi_\pi^N(\k)  . 
\end{equation*}
In particular, if the Jacobi  matrix has a right-limit, i.e.~there exists an (infinite) matrix $\L$ such that for all $i,j \in\mathbb{Z}$, 
\begin{equation*}
\lim_{N\to\infty} \J_{N+i, N+j} = \L_{i,j} 
\end{equation*}
then
\begin{equation}  \label{cumulant_8}
\lim_{N\to\infty}\Cu^n_{K^N_V}[Q]
=  \sum_{m =1}^{M} \sum_{\pi \in \widetilde{\Gamma}_{m}^n } \prod_{i=0}^{n-1}
Q(\L)_{\pi(i)\pi(i+1)}  \sum_{\k \vdash n}\Upsilon_0(\k)\Phi_\pi^0(\k) , 
\end{equation}
where $\widetilde{\mathcal{G}}$ denotes the   adjacency graph of the matrix $Q(\L)$ and
$$
\widetilde{\Gamma}_{m}^n = \big\{ \text{paths }\pi  \text{ on the graph }\widetilde{\mathcal{G}}\text{ of length }n \text{ such that }\pi(0)=\pi(n)=-m  \big\} . 
$$
The condition \eqref{recurrence_limit} implies that the right-limit of the Jacobi matrix is a tridiagonal matrix $\L$ such that $\L_{jj}= 0$ and $\L_{j, j\pm 1} =1/2$   for all $j \in \mathbb{Z}$ and it  was proved in \cite{Lambert_b}, see also~\cite{BD17}, that in this case:
\begin{equation} \label{cumulant_17}
\lim_{N\to\infty} \Cu^n_{K^N_V}[Q]  = \begin{cases}
\displaystyle 
  \sum_{k=1}^{\operatorname{deg} Q}  k \bigg( \int_{-1}^1 Q(x) T_k(x) \frac{dx}{\pi\sqrt{1-x^2}}  \bigg)^2  &\text{if } n=2 \\
 0 &\text{else} 
\end{cases} .
\end{equation}

\medskip

The combinatorial method used in the previous section is well-suited to investigate the global fluctuations of 1-dimensional log-gas, but it is difficult to implement in the mesoscopic regime since we cannot use polynomials as test functions. 
So, to describe the transition for mesoscopic linear statistics and to prove theorem~\ref{thm:crossover_4}, we rely on the asymptotics of the correlation kernel $K^N_V$ from \cite{Deift99+}  and the method from \cite{Lambert_a} that we review below.
Recall that $\varrho_V$ is the equilibrium density of  the gas and define {\it the integrated density of states}:
\begin{equation} \label{IDS}
F_V(x) = \int_0^x  \varrho_V(s) ds . 
 \end{equation}

Let us fix $x_0\in\mathscr{I}_V$, $0<\alpha<1$, and set
\begin{equation} \label{kernel_meso}
\widetilde{K}^N_{V,x_0}(x, y) = \frac{1}{N^\alpha} K_{V}^N\left(x_0+\frac{x}{N^{\alpha}} ,x_0 +\frac{y}{N^{\alpha}}\right)   . 
\end{equation}
Based on the results of \cite{Deift99+},  we have for any $\alpha \in (0,1]$,
\begin{equation} \label{sine_1}
\widetilde{K}^N_{V,x_0}(x, y)   
= \frac{\sin \left[ \pi  N \big(( F_V(x_0+ x N^{-\alpha})-F_V(x_0+ y N^{-\alpha})\big)\right]}{\pi(x-y)}+ \underset{N\to\infty}{O}\left(N^{-\alpha}\right)  ,
\end{equation}
uniformly for all $x, y$ in compact subsets of $\R$; c.f.~\cite[Poposition~3.5]{Lambert_a}. 
The main idea of the method of  \cite{Lambert_a} is to compare the kernel \eqref{sine_1} to the sine-kernel and use the results from Soshnikov \cite{Soshnikov00} for the cumulants of linear statistics of the Sine process. 
We define the sine-kernel with density $\rho>0$ on $\R$ by
\begin{equation} \label{sine}
K^{\sin}_\rho(x,y) = \frac{\sin[\pi \rho(x-y)]}{\pi(x-y)} .
\end{equation}
We see by taking $\alpha=1$ in formula \eqref{sine_1} that  the Sine process with correlation kernel \eqref{sine} describes the local limit in the bulk of the 1-dimensional log-gases. 
In the mesoscopic regime, it was proved in \cite{Lambert_a} that, up to a change of variable, it is possible to replace the kernel $\widetilde{K}^N_{V,x_0}$ by an appropriate sine-kernel using the asymptotics \eqref{sine_1} in the cumulant formulae. Namely, for any $n\ge 2$, 
\begin{equation} \label{sine_2}
\Cu^n_{K^N_V}[f_N] \simeq  \Cu^n_{K^{\sin}_{\eta_N(x_0) }}[ f\circ\zeta_N]
\end{equation}
where  
\begin{equation} \label{zeta_N}
\zeta_N(x) =  N^\alpha \left\{ G_V\left( F_V(x_0) + \varrho_V(x_0) \frac{x}{N^\alpha} \right)-x_0\right\}  . 
\end{equation}
Here, the function $G_V$ denotes the inverse of the integrated density of sates $F_V$, \eqref{IDS}. \
By the inverse function theorem, it exists in a neighborhood of any point $F_V(x_0)$ when $x_0\in\mathscr{I}_V$ and the map  $\zeta_N$ is well-defined  on any compact subset of $\R$ as long as the  parameter $N$ is sufficiently large.
Then, using  Soshnikov's main combinatorial lemma, it was proved in \cite{Lambert_a} that
\begin{equation} \label{cumulant_16}
\lim_{N\to\infty} \Cu^n_{K^N_V}[f_N]   = \begin{cases}  \displaystyle \int_{0}^\infty u \big| \hat{f}(u)  \big|^2 du &\text{if } n=2  \\ 
0  &\text{if } n \ge 3 \end{cases}.
\end{equation}

\subsection{Proof of theorem~\ref{thm:crossover_1} -- The global regime} \label{sect:glob1}

In this section, we modify the strategy described above in order to  deduce theorem~\ref{thm:crossover_1} from our general  theorem~\ref{thm:main}. 
The first step is to verify the assumption \eqref{ass1}. 
By~\cite[Theorem~11.1.2]{PS11}, the expected density of states satisfies for all $x\in\R$,
$$
u^N_V(x) \le e^{- 2N \{ V(x) - \log(1+|x|) - C\} } ,
$$
where $C$ is a constant which depends only on the potential $V$.
This implies that the law of large numbers \eqref{LLN} can be extended to all continuous functions with polynomial growth. 
Moreover, by \eqref{cumulant_17}, we obtain the asymptotics \eqref{CLT0} for the cumulants $ \Cu^n_{K_V^N}[Q]$ of the linear statistic $\Xi(Q)$.  
Then, it only remains to verify that the  estimates \eqref{ass2} hold for any polynomial test function. 
By~\eqref{normalization}, since  $\sum_{\k \vdash n} \Upsilon_m(\k) = 0$,  
the very same computation leading to \eqref{cumulant_8} shows that for any integers $n,m\ge 1$
\begin{align} \notag
\lim_{N\to\infty} \underset{x_{0} =x_n}{\int_{\X^n}}     \Upsilon^n_m[Q](x)  \prod_{1\le j\le n} K^N_V(x_j, x_{j-1}) d^nx 
=  \sum_{m =1}^{M} \sum_{\pi \in \widetilde{\Gamma}_{m}^n } \prod_{i=0}^{n-1}
Q(\L)_{\pi(i)\pi(i+1)}  \sum_{\k \vdash n}\Upsilon_m(\k)\Phi_\pi^0(\k) . 
\end{align}
Since the (infinite) matrix $\L$ is bounded with $\|\L\| \le 1/2$, we obtain the estimates \eqref{ass2} with an error which is $O(1)$. This completes the proof of theorem~\ref{thm:crossover_1}.

\begin{remark}[Generalizations of theorem~\ref{thm:crossover_1}] \label{rk:biorthogonal}
Note that we have formulated theorem~\ref{thm:crossover_1} for a log-gas at inverse temperature $\beta=2$, but the previous proof can be  generalized to other one-dimensional biorthogonal ensemble with a correlation kernel of the form:
  \begin{equation*} 
K^N(z,w) = \sum_{k=0}^{N-1} \varphi_k^N(z) \varpi_k^N(w) .  
\end{equation*} 
The appropriate assumptions are that there exists an equilibrium density and a law of large numbers holds for all polynomials,  the family $\{\varphi_k^N \}_{k=0}^\infty $ satisfies a $q$-term recurrence relation for all $N\in\N$  and the corresponding recurrence matrix $\J$  has a right-limit~$\L$ as $N\to\infty$.  
 This applies to other orthogonal polynomial ensembles, such as the discrete point processes coming from domino tilings of hexagons, as well as some  non-symmetric biorthogonal  ensembles  such as the Muttalib-Borodin ensembles, square singular values of product of complex Ginibre matrices or certain two-matrix models,  see \cite{BD17, Lambert_b} for more details.
 Moreover, we only require that the right-limit $\L$ exists but it need not be a Toeplitz matrix. 
 Then, we obtain a crossover from a non-Gaussian process (described by $\L$ in the regime where $T_N\to 0$) to a Poisson process (when $T_N\to \infty$). 
For instance, such a transition arises when considering linear statistics of the log-gases in the multi-cut regime, \cite{Lambert_b}    
%Note also that  formula    \eqref{cumulant_16} still  holds in the case where the right-limit $\L$ exists but is not Toeplitz (i.e.~constant along its main diagonals). Then, we obtain a crossover from a non-Gaussian process (described by $\L$ in the regime where $T_N\to 0$) to a Poisson process (when $T_N\to \infty$). 
%For instance, such a transition arises when considering linear statistics of the log-gases in the multi-cut regime.    
 \end{remark}

\begin{remark}\label{rk:one_cut}
It was proved in \cite[section~5]{Johansson98} that when $V$ is a convex polynomial, then $\S_V=(-1,1)$ and the conditions~\eqref{recurrence_limit} are satisfied. In fact, Johansson's argument shows that these conditions are also necessary to have a CLT for polynomial test functions.   Therefore, it is an interesting question to know whether \eqref{recurrence_limit} and the one-cut condition $\S_V=(-1,1)$ are equivalent.
\end{remark}

\subsection{Proof of theorem~\ref{thm:crossover_4} -- The mesoscopic regime} \label{sect:meso1}

Let us fix $x_0\in\mathscr{I}_V$, $0<\alpha<1$ and let $f_N= f\big(N^\alpha( \cdot -x_0) \big)$ where $f\in \mathcal{C}^2_c(\R)$. 
First, observe that by a change of variable, we have for all $n\ge 1$, 
    \begin{equation*} 
 \int_\R f_N(x)^n  K_V^N(x,x) dx =   \int_\R f(x)^n   \widetilde{K}^N_{V,x_0}(x,x) dx ,
  \end{equation*}
where $\widetilde{K}^N_{V,x_0}$  is given by \eqref{kernel_meso}. 
Using the asymptotics \eqref{sine_1}, we have that 
$$\widetilde{K}^N_{V,x_0}(x, x)  = N^{1-\alpha} \varrho_V(x_0)  + \underset{N\to\infty}{O}\left(N^{-\alpha}\right)  $$
uniformly for all $x\in\supp(f)$. If $M_N = N^{1-\alpha} \varrho_V(x_0) $, this implies that for all $n\in\N$, 
\begin{equation*}
\frac{1}{M_N} \int_\R f_N(x)^n  K_V^N(x,x) dx
\simeq  \int_\R  f(x)^{n} dx .
\end{equation*}
Thus, we obtain the condition \eqref{ass1} with $d\eta =dx$. 
Moreover, the assumption \eqref{CLT0} is given by \eqref{cumulant_16}.  Then it just remains to prove the estimates \eqref{ass2}. 
Let us observe that by a change of variables, we also have for all $n, m \ge 1$, 
\begin{equation} \label{cumulant_13}
\underset{u_{0} =u_n}{\int_{\R^n}}   \Upsilon^n_m[f_N](u)  \prod_{1\le j\le n} K_V^N(u_j, u_{j-1}) d\mu^n(u) 
= \underset{u_{0} =u_n}{\int_{\R^n}}   \Upsilon^n_m[f](u)  \prod_{1\le j\le n} \widetilde{K}^N_{V,x_0}(u_j, u_{j-1}) d\mu^n(u) . 
\end{equation}
Exactly like the proof of  \eqref{sine_2} -- which corresponds to the case $m=0$ by formula \eqref{cumulant_0} --  we deduce from the proof of \cite[proposition~2.2]{Lambert_a} that  for any $m \ge 1$, 
\begin{equation}\label{cumulant_10}
 \underset{x_{0} =x_n}{\int_{\R^n}}     \Upsilon^n_m[f](x)  \prod_{j\le n} \widetilde{K}^N_{V,x_0}(x_j, x_{j-1}) d^nx \
 \simeq \underset{x_{0} =x_n}{\int_{\X^n}}  \hspace{-.2cm}   \Upsilon^n_m[h_N](x)  \prod_{1\le j\le n} K^{\sin}_{\eta_N}(x_j, x_{j-1}) d^nx
\end{equation}
where $h_N = f\circ \zeta_N$ and $\zeta_N$ is given by \eqref{zeta_N}. 
To finish the proof, we need the following estimates.

\begin{proposition} \label{thm:g}
%Let  $x_0\in \mathscr{I}_V$, $0<\alpha<1$, $f \in C^2_0(\R)$.
Suppose that $\supp(f) \subset (-L, L)$.  There exists 
$N_0 >0$ such that for all $N \ge N_0$,
 the functions $h_N= f\circ \zeta_N$ are well-defined on $\R$,  $ h_N\in  C^2_c([-L,L])$ and for all $u\in\R$, 
\begin{equation} \label{g_estimate} 
\big| \widehat{h_N}(u) \big| \le \| f \|_{\mathcal{C}^2(\R)} \frac{C }{1+ |u|^2} .
\end{equation} 
\end{proposition}

\proof When the potential $V$ is analytic, the {\it bulk} $\mathscr{I}_V$ consists of finitely many open intervals and the equilibrium density $\varrho_V$ is smooth on $\mathscr{I}_V$. Since $x_0\in\mathscr{I}_V$, by formula \eqref{zeta_N}, the function $\zeta_N$ is increasing and smooth on the interval $[-L,L]$ with
$$  \zeta_N''(x) =   \varrho_V(x_0)^2 G_V''\big(F_V(x_0) + \varrho_V(x_0) xN^{-\alpha}\big ) N^{-\alpha}. 
$$
Moreover, since $\zeta_N(0)=0$ and  $\zeta'_N(0) =G_V'\big(F_V(x_0)\big) \varrho_V(x_0) =1$, this implies that 
$$  \zeta_N(x) = x + \O(N^{-\alpha}) $$  
uniformly for all $x\in[-L,L]$.
Since the open interval $(-L, L)$ contains the support of the test function $f$, this estimate shows that when the parameter $N$ is large, we can define  $h_N(x) = f\big( \zeta_N(x) \big)$ for all $x\in[-L, L]$ and extend it by 0 on $\R \backslash [-L,L]$.  
Then  $h_N \in  C^2_0(\R)$ and 
\begin{equation} \label{h_N''}
h_N''(x) =  \zeta_N''(x) f'(\zeta_N(x)) +  \zeta_N'(x)^2f''(\zeta_N(x)) 
\end{equation}
for all $x\in[-L,L]$. Moreover, we can use the estimate
\begin{equation} \label{Fourier_estimate} 
\big| \widehat{h_N}(u) \big| \le C \frac{\| h_N \|_\infty + \| h_N''\|_\infty}{1+ |u|^2} 
\end{equation}
to get the upper-bound \eqref{g_estimate}. Plainly $\| h_N \|_\infty  \le \| f\|_\infty$ and it is easy to  deduce from formula \eqref{h_N''} that
$$ \| h_N'' \|_\infty \le \| f' \|_\infty + C \|f''\|_\infty . $$  
\qed

To compute the limit of the RHS of \eqref{cumulant_10}, we also need the following asymptotics which come from the proof of Lemma~1 in Soshnikov's paper \cite{Soshnikov00} on linear statistics of the CUE and Sine process.

\begin{lemma} \label{thm:sine}
Let $n\ge 2$ and let $\eta_N>0$ such that $\eta_N\nearrow \infty$ as $\to+\infty$.   Suppose that $h_N$ is a sequence of integrable functions such that 
\begin{equation}\label{assumption_2}
\lim_{N\to\infty} \hspace{-.3cm}
\underset{\begin{subarray}{c} u_1+\cdots +u_{n}=0  \\ |u_1|+\cdots +|u_n| > \eta_N \end{subarray}}{\int_{\R^{n-1}}} \hspace{-.4cm} \big| \widehat{h_N}(u_1)\cdots   \widehat{h_N}(u_n) \big| |u_1|  d^{n-1}\mathrm{u} =0 . 
\end{equation}
Then,  for any map $\Upsilon: \mho \to\R$ such that $\sum_{\k \vdash n} \Upsilon(\k) =0$, we have
\begin{equation}\label{cumulant_11}
 \hspace{-.5cm}
\underset{x_{0} =x_n}{\int_{\X^n}}  \hspace{-.2cm}   \Upsilon^n[h_N](x)  \prod_{1\le j\le n} K^{\sin}_{\eta_N}(x_j, x_{j-1}) d^nx\
\simeq -
 \hspace{-.3cm}\underset{u_1+\cdots +u_n=0}{\int_{\R^{n-1}}} \hspace{-.3cm}  \Re\bigg\{ \prod_{j=1}^n  \widehat{h_N}(u_j)  \bigg\}\sum_{\k \vdash n}\Upsilon(\k)\Psi_u(\k)  d^{n-1}u ,
\end{equation}
where for  any $u \in \R^n$ and for any composition $\k\vdash n$, 
$$ \Psi_u(\k) =  2\max_{1\le j <\ell(\k)}\{0, u_1+\cdots + u_{k_1+\cdots +k_j} \} .
$$

\end{lemma}

\begin{proof}
Based on the formula
$$
K^{\sin}_{\eta_N}(x,y) = \int_\R \1_{\{ |u|< \eta_N/2\} } e^{2\pi i u (x-y)}  du , 
$$
we obtain 
\begin{align*}
\mathscr{T}_N&:=\underset{x_{0} =x_n}{\int_{\X^n}}  \hspace{-.2cm}   \Upsilon^n[h_N](x)  \prod_{1\le j\le n} K^{\sin}_{\eta_N}(x_j, x_{j-1}) d^nx \\
& = 
 \hspace{-.2cm}\underset{u_1+\cdots +u_n=0 \hspace{.1cm}}{\int_{\R^{n-1}}}  \prod_{j=1}^n  \widehat{h_N}(u_j)  \sum_{\k \vdash n}\Upsilon(\k) \max\big\{0, \eta_N - \Psi_u(\k)/2 - \Psi_{-u}(\k)/2 \big\} d^{n-1}u .
\end{align*}
Then, the condition $\sum_{\k \vdash n} \Upsilon(\k) =0$ implies that 
\begin{align*}
&\bigg| \mathscr{T}_N +
 \hspace{-.2cm}\underset{u_1+\cdots +u_n=0 \hspace{.1cm}}{\int_{\R^{n-1}}}  \prod_{j=1}^n  \widehat{h_N}(u_j)  \sum_{\k \vdash n}\Upsilon(\k) \frac{ \Psi_u(\k) + \Psi_{-u}(\k)}{2} d^{n-1}u \bigg| \\
 & \le 
 \hspace{-.2cm}\underset{u_1+\cdots +u_n=0 \hspace{.1cm}}{\int_{\R^{n-1}}}
   \bigg| \prod_{j=1}^n  \widehat{h_N}(u_j) \bigg|  \sum_{\k \vdash n} |\Upsilon(\k) |  \Psi_u(\k)    \1_{ \{\Psi_u(\k) + \Psi_{-u}(\k) \ge 2 \eta_N \} }d^{n-1}u 
\end{align*}

Since $ \big| \Psi_u(\k)/2 \big| \le |u_1|+\cdots +|u_n|$ for any $\k\vdash n$, the condition \eqref{assumption_2} is sufficient to obtain the asymptotics \eqref{cumulant_11}. 
\end{proof}

\medskip

We are now ready to complete the proof of theorem~\ref{thm:crossover_4}.
Using the estimate \eqref{g_estimate}, we see that when the parameter $N$ is sufficiently large,  there exists a constant $C>0$ so that  
\begin{equation} \label{cumulant_12}
\begin{aligned} 
&\Bigg| \underset{u_1+\cdots +u_n=0}{\int_{\R^{n-1}}} \hspace{-.3cm}  \Re\bigg\{ \prod_{j=1}^n  \widehat{h_N}(u_j)  \bigg\}\sum_{\k \vdash n}\Upsilon_m(\k)\Psi_u(\k)  d^{n-1}u  \Bigg| \\
& \hspace{2cm}
\le C \hspace{-.2cm}  \underset{u_1+\cdots +u_n=0}{\int_{\R^{n-1}}}  \frac{|u_1|+\cdots +|u_n|}{(1+ |u_1|^2) \cdots (1+|u_n|^2)} d^{n-1}u .
\end{aligned} 
\end{equation}
A similar upper-bound shows that the sequence $h_N = f\circ \zeta_N$ satisfies the condition \eqref{assumption_2} of lemma~\ref{thm:sine}. Hence, combining the asymptotics \eqref{cumulant_10} and \eqref{cumulant_11}, we obtain
\[
 \underset{x_{0} =x_n}{\int_{\R^n}}     \Upsilon^n_m[f](x)  \prod_{j\le n} \widetilde{K}^N_{V,x_0}(x_j, x_{j-1}) d^nx  \simeq  -
 \hspace{-.3cm}\underset{u_1+\cdots +u_n=0}{\int_{\R^{n-1}}} \hspace{-.3cm}  \Re\bigg\{ \prod_{j=1}^n  \widehat{h_N}(u_j)  \bigg\}\sum_{\k \vdash n}\Upsilon_m(\k)\Psi_u(\k)  d^{n-1}u .
\]
By \eqref{cumulant_12},  we see that the previous integral is uniformly bounded by a constant which depends only on the test function $f$ and $n,m \in\N$. 
Hence, by \eqref{cumulant_13}, we obtain the estimates \eqref{ass2} with an error which is $O(1)$ and we can apply theorem~\ref{thm:main}. This completes the proof of  theorem~\ref{thm:crossover_4}.

\begin{appendices}

\appendix\section{Incomplete determinantal processes} \label{A:kernel}

In this appendix, we review some background material on determinantal processes and we give an alternative proof of the results of Bohigas and Pato, \cite{BP06}, that the incomplete process is determinantal with correlation kernel $p_NK^N_V(z,w)$. 
There are many excellent surveys about  determinantal processes and we refer to 
\cite{Soshnikov, Konig, Johansson06, Borodin11} for further details and an overview of the  main examples. 

\medskip

Let $\X$ be a complete separable metric space equipped with a Radon measure~$\mu$. 
The configuration space $\mathscr{Q}(\X)$ is the set of integer-valued locally finite Borel measure on $\X$ equipped with the topology which is generated by the maps 
$$\mathscr{Q}(\X) \ni \Xi \mapsto \Xi_A := \int_A d\Xi    $$ for every Borel set $A\subseteq \X$.
A point process $\mathbb{P}$ is a Borel probability measure on  the configuration  space  $\mathscr{Q}(\X)$. 
A point process can be characterized by its correlation functions $(\rho_n)_{n=1}^\infty$. If it exists, $\rho_n$ is a symmetric function on $\X^n$ which satisfies the identity   
\begin{equation}\label{correlation}
\E\left[  {\Xi_{A_1} \choose k_1} \cdots   {\Xi_{A_\ell} \choose k_\ell} \right] = \frac{1}{k_1!\cdots k_\ell!} \int_{A_1^{k_1} \times \cdots \times A_\ell^{k_\ell}} \rho_n(x_1,\dots , x_n) d\mu^n(x)
\end{equation}
 for all compositions $\k \vdash n$  and for all disjoint Borel sets $A_1,\dots, A_\ell \subseteq \X$. 
A point process is called {\bf determinantal} if all its correlation functions exist and are given by
\begin{equation} \label{det}
\rho_n(x_1,\dots , x_n) = \det_{n\times n}\big[ K(x_i, x_j) \big] . 
\end{equation}
The function $K: \X \times \X \to\C$ is called the {\bf correlation kernel}. 
For instance,  given random points $(\lambda_1,\dots, \lambda_N)$ with the joint density $G_N(x)=e^{-\beta \mathscr{H}^N_V(x)}/Z^N_V$ on $\X^N$, see \eqref{log_gas}, the random measure $ \sum_{k=1}^N \delta_{\lambda_k}$ defines a point process  and its correlation functions  satisfy $\rho_N = N! G_N$ and 
\begin{equation}\label{correlation_function}
\rho_k(x_1, \dots, x_k) = \frac{N!}{(N-k)!}   \int  G_N(x_1, \dots, x_N) dx_{k+1}\cdots dx_N
\end{equation}
for all $k<N$. When $\beta=2$, it is easy to verify that 
$$ G_N(x) = \frac{1}{N!} \det_{N\times N} \big[K^N_V(x_i, x_j) \big]$$ 
and that, by formula \eqref{correlation_function},  the process $\Xi$ is determinantal with the correlation kernel 
$K^N_V$ given by \eqref{kernel}. In general, if  $K$  is a continuous, Hermitian symmetric, function which satisfies property~\eqref{projection}, then there exists a determinantal process on $\X$ with correlation kernel $K$; c.f.~\cite[Theorem~3]{Soshnikov}. 

\medskip

%In this general setting, let $\Xi and  

Let $q\in (0,1)$ and $p=1-q$.
Recall that $X$ is a Binomial random variable with parameter $p$ and $N\in\N$ if for all $k \in \N_0$,  
\begin{equation} \label{binomial_moment}
 \E\left[{X \choose k} \right]   = p^k {N \choose k} . 
 \end{equation}
By convention, ${N \choose k}  = 0$ if $k>N$. 
Let $\Xi$ be a random (point) configuration and let $\widehat{\Xi}$ be the configuration obtained after performing a Bernoulli percolation on  $\Xi$. 
By construction,  for any disjoint Borel set $A\subseteq \X $, the conditional distribution of   the random variable $\widehat{\Xi}_A$ given $\Xi$  has a Binomial distribution with parameters $p$ and $\Xi_A$ and it is statistically independent of $\Xi_B$ for any Borel set $B$ disjoint of $A$. 
By formula \eqref{binomial_moment}, this implies that
\begin{align*}
\E\left[  {\widehat\Xi_{A_1} \choose k_1} \cdots   {\widehat\Xi_{A_\ell} \choose k_\ell} \bigg|  \Xi_A=N_1, \dots, \Xi_A=N_\ell \right]
= p^n  {N_1 \choose k_1} \cdots   {N_\ell \choose k_\ell} 
\end{align*}
 for any composition $\k \vdash n$  and for all disjoint Borel sets $A_1,\dots, A_\ell \subseteq \X$.
Hence, we obtain
\begin{align*}
\E\left[  {\widehat\Xi_{A_1} \choose k_1} \cdots   {\widehat\Xi_{A_\ell} \choose k_\ell} \right]
= p^n \E\left[  {\Xi_{A_1} \choose k_1} \cdots   {\Xi_{A_\ell} \choose k_\ell} \right] 
\end{align*}
so that,  by formula \eqref{correlation}, the correlation functions  of {\it the incomplete process} $\widehat{\Xi}$ are given by $p^n \rho_n(x_1,\dots , x_n)$ for all $n \ge 1$. 
In particular, we deduce from formula \eqref{det} that, if  $\Xi$ is a determinantal process with a correlation kernel $K$,  then the point process $\widehat{\Xi}$  is also determinantal with kernel $pK$.  

%  In particular, the  space  $\mathscr{Q}(\X)$ is equipped with its Borel sigma-algebra a natural notion of weak convergence. 

\section{Off-diagonal decay of the correlation kernel $K^N_V$ in dimension 2} \label{A:damping}

In this section, we review some classical estimates for the correlation kernel~\eqref{kernel} which  have been used in~\cite{AHM11} to prove the CLT \eqref{CLT_0}.  
Then, we prove lemma~\ref{thm:localization} and an analogous result for the cumulants of the \Ginibre process. For completeness, we also give the proof of lemma~\ref{thm:Ginibre_approximation}.  
We will use the formulation of section~5  in \cite{AHM11} but the estimates \eqref{exp_decay} and \eqref{tail_bound} go back to the papers~\cite{Berman09} and~\cite{AHM10}.  Suppose that the potential $V:\C \to\R$ is real-analytic and satisfies the condition \eqref{potential}. 
Then, there a function $\phi_V : \C \to \R^+$ such that $\phi_V(z) \ge \nu \log |z|^2$ as $|z| \to\infty$ and some constants $C, c, \delta>0$ such that
\begin{equation}\label{exp_decay}
\big| K_V^N(z,w) \big| \le C N e^{- c\sqrt{N} (|z-w|\wedge \delta )}  ,
\end{equation}
and 
\begin{equation}\label{tail_bound}
K_V^N(z,z) \le C N e^{-N  \phi_V(z)} , 
\end{equation}
for all $w \in \S$ and $z\in \C$. \\ % It was observed in~\cite{AHM11} that these estimates implies lemma~\ref{thm:localization}. 

{\it Proof of lemma~\ref{thm:localization}.} We use  the convention  $z_0=z_{n+1}$.  Since the kernel $K_V^N$ is reproducing, by the Cauchy-Schwartz  inequality, 
$$ \big| K_V^N (z,w) \big| \le \sqrt{ K_V^N (z,z) K_V^N (w,w) } , $$
for all $z, w\in\C$, so that
$$
\bigg| \prod_{j=0}^{n}K_V^N(z_j,z_{j+1}) \bigg|  \le \big| K_V^N(z_0,z_{1})  \big|\sqrt{K_V^N(z_0,z_{0})K_V^N(z_1,z_{1})} \prod_{j=2}^{n}K_V^N(z_j,z_j) .
$$

Since $\displaystyle \int_\C K_V^N(z,z) d\A(z) = N$ and $K_V^N(z,z)  \le C N $ (see the estimate \eqref{tail_bound}), we obtain 
$$
 \underset{|z_1-z_{n+1}|> \epsilon_N}{\int_{\C^n\times \S}}  \bigg|\prod_{j=0}^n K^N_V(z_{j}, z_{j+1})\bigg|\   d\A^{n+1}(\z)  \le  C N^n  \underset{|z_1-z_0|> \epsilon_N}{\int_{\S\times \C}} \hspace{-.2cm}  \big| K_V^N(z_0 ,z_1) \big| d\A(z_0)   d\A(z_1)  . 
$$
Then, it easy to check that the estimates  \eqref{exp_decay} and \eqref{tail_bound} imply that 
$$
 \underset{|z_1-z_0|> \epsilon_N}{\int_{\S \times \C}}   \big| K_V^N(z_0 ,z_1) \big| d\A(z_0)   d\A(z_1)  
\le C N  e^{-c \sqrt{N} \epsilon_N }  . 
$$
Hence, if $\epsilon_N = \kappa N^{-1/2} \log N $ and $\kappa \ge (n+1)/ c$, we obtain 
\begin{equation}\label{estimate_6}
 \underset{|z_1-z_{n+1}|> \epsilon_N}{\int_{\C^n \times \S}} \bigg| \prod_{j=0}^n K^N_V(z_{j}, z_{j+1})\bigg|\   d\A^{n+1}(\z)  = \O(N^{-1}) .
\end{equation}

Moreover, since  $\sup\big\{ |F_N(z_0, \z)| : \z \in \C^{n} , N  \ge N_0 \big\}  \le C \1_{z_0 \in\S}$ by  \eqref{estimate_0}, the estimate \eqref{estimate_6} implies that 
\begin{align} \label{estimate_4}
&\underset{z_0 = z_{n+1}}{\int_{\C^{n+1}}}\hspace{-.2cm} F_N(\z) \prod_{j=0}^n K^N_V(z_{j}, z_{j+1}) d\A^{n+1}(\z)  \\
&\notag\hspace{1cm}
 =  \underset{|z_1-z_{n+1}| \le \epsilon_N}{\int_{\C^n \times \S}} \hspace{-.2cm}  F_N(\z)   \prod_{j=0}^n K^N_V(z_{j}, z_{j+1}) \   d\A^{n+1}(\z)  + \O(N^{-1}) .
\end{align}

Now, we can proceed by induction to get formula \eqref{localization}. If
 $\mathscr{C}_N = \{\z \in \C^{n+1} :  z_{n+1} \in \S, |z_1-z_{n+1}| \le \epsilon_N \} $, the next step is to show that 
\begin{equation} \label{estimate_2}
 \underset{|z_2-z_1|> \epsilon_N}{\int_{\mathscr{C}_N}} \bigg| \prod_{j=0}^n K^N_V(z_{j}, z_{j+1})\bigg|\   d\A^{n+1}(\z)  = \underset{N\to\infty}{O}(N^{-1}) .
\end{equation}

Since the set $\S_V$ is open, there exists a compact set $\S' \subset \S_V$ such that $\S \subset \S'$ and  $\mathscr{C}_N \subset  \{\z \in \C^{n+1} :  z_{n+1} , z_1 \in \S' \}$  when the parameter $N$ is sufficiently large. Then, as before, we obtain
$$
 \underset{|z_2-z_{1}|> \epsilon_N}{\int_{\mathscr{C}_N}}  \bigg|\prod_{j=0}^n K^N_V(z_{j}, z_{j+1})\bigg|\   d\A^{n+1}(\z)  \le  C N^n  \underset{|z_2-z_1|> \epsilon_N}{\int_{\S'\times \C}} \hspace{-.2cm}  \big| K_V^N(z_1 ,z_2) \big| d\A(z_0)   d\A(z_1)  , 
$$
and  formula \eqref{estimate_2} also follows directly from the estimate \eqref{estimate_6}.
Hence, by formula \eqref{estimate_4}, this implies that
\begin{align*}
&\underset{z_0 = z_{n+1}}{\int_{\C^{n+1}}}\hspace{-.2cm} F_N(\z) \prod_{j=0}^n K^N_V(z_{j}, z_{j+1}) d\A^{n+1}(\z)  \\
&\notag\hspace{1cm}
 =  \underset{\begin{subarray}{c} |z_1-z_{n+1}| \le \epsilon_N  \\  |z_2-z_{1}| \le \epsilon_N  \end{subarray}} {\int_{\C^n \times \S}} \hspace{-.2cm}  F_N(\z)   \prod_{j=0}^n K^N_V(z_{j}, z_{j+1}) \   d\A^{n+1}(\z)  + \O(N^{-1}) .
\end{align*}
If we repeat this argument, we obtain \eqref{localization}. \qed

\medskip
  
 \begin{lemma} \label{thm:delocalization}
Let $n \in\N$, $w_0= w_{n+1} =0$, and let $H(\w)$ be a polynomial of degree at most 2 in the variables $w_1,\dots, w_n, \bar{w_1},\dots, \bar{w_n}$.
 For any  sequence $\delta_N \ge k  \rho_N^{-1/2} \sqrt{\log\rho_N}$ with $k>0$, we have
\begin{align*}
& \int_{\C^n} H(\w) \prod_{j=0}^nK^\infty_{\rho_N}(w_{j}, w_{j+1})\  d\A^n(\w)  \\
  &\hspace{1cm}=  \int_{\mathscr{A}(0; \delta_N)} \hspace{-.2cm} H(\w) \prod_{j=0}^nK^\infty_{\rho_N}(w_{j}, w_{j+1})\   d\A^n(\w)  
+ \O(\rho_N^{1-k^2/2})
\end{align*}
where the set $\mathscr{A}(0; \delta_N)$ is given by formula \eqref{A}. 
 \end{lemma}
 
 \proof We will first  show that
\begin{equation} \label{estimate_1}
  \underset{|w_1| > \delta_N}{\int_{\C^n}}   H(\w) \prod_{j=0}^nK^\infty_{\rho_N}(w_{j}, w_{j+1}) \  d\A^n(\w)
   = \underset{N\to\infty}{O}\big( \rho_N e^{- \rho_N \delta_N^2/2}    \big) .
\end{equation}
 First, notice that if $H =1$, since $w_0=0$, by formula~\eqref{estimate_5}, we have 
  $$ 
  \underset{|w_1| > \delta_N}{\int_{\C^n}} \bigg| \prod_{j=0}^nK^\infty_{\rho_N}(w_{j}, w_{j+1})\bigg| \  d\A^n(\w) 
   \le  e^{- \rho_N \delta_N^2/2}    \rho_N^{n+1}
   \underbrace{\int_{\C^n} \prod_{j=1}^n e^{-\rho_N |v_j|^2/2}  d\A^n(\mathrm{v})}_{\displaystyle = \rho_N^{-n}} .    $$
 
In general, there exists a constant $C>0$ which only depends  on the polynomial $H$ so that 
$$
\big| H(\w) \big|  \le  C \big\{ 1+ |w_1|^2 +\cdots +|w_n|^2 \big\} 
$$
or 
  \begin{equation} \label{estimate_3}
\big| H(\w) \big|  \le  C \big\{ 1+ |w_2-w_1|^2 +\cdots +|w_{n-1} -w_n|^2 + |w_n|^2 \big\} . 
\end{equation}
Since, for any $k=1,\dots, n$,
\begin{align*}
&    \underset{|w_1| > \delta_N}{\int_{\C^n}}|w_k-w_{k+1}|^2 \bigg|  \prod_{j=0}^nK^\infty_{\rho_N}(w_{j}, w_{j+1})  \bigg| d\A^n(\w) \\
& \hspace{2cm}   \le  e^{- \rho_N \delta_N^2/2} 
   \underbrace{    \rho_N^{n+1} \int_{\C^n}  |v_k|^2 \prod_{j=1}^n e^{-\rho_N | v_{j} |^2/2}  d\A^n(\mathrm{v})}_{\displaystyle = 1} ,
   \end{align*}
the estimate \eqref{estimate_1} follows  directly from  \eqref{estimate_3} and the leading contribution comes from the constant term. If we use the estimate 
\begin{equation*} \big| H(\w) \big|  
\le  C \bigg\{ 1+ |w_1|^2 + \sum_{\begin{subarray}{c} j=1 \\ j \neq k \end{subarray}}^n |w_{j+1} -w_j|^2 +|w_n|^2 \bigg\} 
\end{equation*}
instead, the same argument shows that for any $k= 1,\dots n$,   
  \begin{equation*}
  \int\limits_{|w_k -w_{k+1}| > \delta_N } \hspace{-.3cm} H(\w) \prod_{j=0}^nK^\infty_{\rho_N}(w_{j}, w_{j+1}) \  d\A^n(\w)
   = \underset{N\to\infty}{O}\big( \rho_N e^{- \rho_N \delta_N^2/2}    \big) .
\end{equation*}
Hence, the lemma follows from applying a union bound and from the choice of the sequence $\delta_N$.  \qed
 
 \medskip

{\it Proof of lemma~\ref{thm:Ginibre_approximation}.}
The map $(z, w) \mapsto \Psi(z,w)$ is bi-holomorphic in a neighborhood of $(x_0, \bar{x_0})$, so there exists $0<\epsilon<1$ such that for all $|u|, |v| \le \epsilon$,
$$
\Psi(x_0 + u, \bar{x_0} + v) =
\sum_{k, j \ge 0} a_{kj} u^k v^j .
$$

By definition, 
$\overline{\Psi(z,w)} = \Psi(\overline{w},\overline{z})$, so that the coefficients of the previous power series are Hermitian-symmetric: $a_{kj} = \bar{a_{jk}}$ for all $k, j \ge 0$. Moreover, by definition,  we have
\begin{equation}\label{density}
a_{11} = \p_z \pb_z V |_{z=x_0} = \Delta V (x_0) = \frac{b_0(x_0, \bar{x_0})}{2} . 
\end{equation}

Let 
$$\h(u) = i \sum_{k>0} \big\{  a_{k0} u^k - a_{0k} \bar{u}^k \big\}   = -2  \Im\bigg\{  \sum_{k>0} a_{k0} u^k   \bigg\} . $$
 Since $V(x_0 + u) = \Psi(x_0 + u , \bar{x_0} + \bar{u} )$, we see that for any $|u|, |v| \le \epsilon$, 
\begin{align*}
& 2\Psi(x_0 + u, \bar{x_0} +\bar{v}) -V(x_0+u)- V(x_0+v)   \\
&\hspace{2cm}= - i\{ \h(u) - \h(v)\} +  a_{11} \big( 2 u \bar{v} - |\bar{u}|^2 -|\bar{v}|^2   \big)    
 + O( \epsilon^3) . 
\end{align*}
By formula \eqref{B_kernel}, this implies that  for any $|u|, |v| \le \epsilon_N = \log (N^\kappa) N^{-1/2}$, 
\begin{equation*}
\frac{B^N(x_0 + u , x_0 +v)e^{i N \h(u)}}{e^{i N \h(v)}}=  N b_0(x_0, \bar{x_0}) e^{N a_{11}  ( 2 u \bar{v} - |\bar{u}|^2 -|\bar{v}|^2)   } \left\{ 1 + \underset{N\to\infty}{O}\big(  (\log N)^{2}\epsilon_N   \big)  \right\} .
\end{equation*}
By formula \eqref{density} and the definition of the  \Ginibre kernel, it completes the proof. \qed
\end{appendices}


\medskip



\begin{thebibliography}{10}

\bibitem{AHM10}
{\sc Y.~Ameur, H.~Hedenmalm, and N.~Makarov}, {\em Berezin transform in
  polynomial {B}ergman spaces}. Comm. Pure Appl. Math.~63, (2010),
  1533--1584.

\bibitem{AHM11}
{\sc Y.~Ameur, H.~Hedenmalm, and N.~Makarov}, {\em Fluctuations of
  eigenvalues of random normal matrices}. Duke Math. J.~159, (2011),
  31--81.

\bibitem{AHM15}
{\sc Y.~Ameur, H.~Hedenmalm, and N.~Makarov}, {\em Random normal matrices and
  {W}ard identities}. Ann. Probab. 43, (2015), 1157--1201.

\bibitem{AGZ}
{\sc G.~W. Anderson, A.~Guionnet, and O.~Zeitouni}, {\em An Introduction to
  Random Matrices}. Cambridge Univ. Press, Cambridge (2010).

\bibitem{BL16}
{\sc F.~Bekerman and A.~Lodhia}, {\em Mesoscopic central limit theorem for
  general $\beta$-ensembles}. 
  Ann. Inst. Henri PoincarŽ Probab. Stat. 54 (2018), no. 4, 1917Ð1938.
%\newblock arXiv:1605.05206

\bibitem{Berggren17}
{\sc T.~Berggren and M.~Duits},
{\em Mesoscopic fluctuations for the thinned Circular Unitary Ensemble}.
Math. Phys. Anal. Geom. 20 (2017), no. 3, Art. 19, 40pp.
%\newblock arXiv:1611.00991

\bibitem{Berman09}
{\sc R.J. Berman}, {\em Bergman kernels for weighted polynomials and weighted
  equilibrium measures of {$\C^n$}}. Indiana U. Math. J. 58, (2009),
  1921--1946.

\bibitem{BP04}
{\sc O.~Bohigas and M.~P. Pato}, {\em Missing levels and correlated spectra}. Phys. Lett. B 595, (2004),  036212,  171--176.

\bibitem{BP06}
{\sc O.~Bohigas and M.~P. Pato}, {\em Randomly incomplete spectra and
  intermediate statistics}. Phys. Rev. E (3) 74, (2006), 036212.

\bibitem{BBNY16}
{\sc R.~Bauerschmidt, P.~Bourgade, M.~Nikula and H.-T. Yau},
{\em The two-dimensional Coulomb plasma: quasi-free approximation and central limit theorem}. \newblock arXiv:1609.08582

\bibitem{Borodin11}
{\sc A.~Borodin}, {\em Determinantal point processes}, in {\em The Oxford handbook of random matrix theory}. Oxford Univ. Press, Oxford, 2011, 231--249.

\bibitem{BG13a}
{\sc G.~Borot and A.~Guionnet}, {\em Asymptotic expansion
  of beta matrix models in the one-cut regime}. Comm. Math. Phys. 317 (2),
  (2013), 447--483.

\bibitem{BG13b}
{\sc G.~Borot and A.~Guionnet}, {\em Asymptotic expansion of beta matrix models
  in the multi-cut regime}.
\newblock arXiv:1303.1045.

\bibitem{BDIK15}
{\sc T.~Bothner, P.~Deift, A.~Its, and I.~Krasovsky}, {\em On the asymptotic
  behavior of a log gas in the bulk scaling limit in the presence of a varying
  external potential {I}}. Comm. Math. Phys. 337, (2015), 1397--1463.

\bibitem{BDIK16}
{\sc T.~Bothner, P.~Deift, A.~Its, and I.~Krasovsky}, {\em On the asymptotic
  behavior of a log gas in the bulk scaling limit in the presence of a varying
  external potential II}.
 In  Large truncated Toeplitz matrices, Toeplitz operators, and related topics, 213--234, Oper. Theory Adv. Appl., 259, BirkhŠuser/Springer, Cham, 2017
%\newblock arXiv:1512.02883 


\bibitem{BD16}
{\sc J.~Breuer and M.~Duits}, {\em Universality of mesoscopic fluctuations for
  orthogonal polynomial ensembles}. Comm. Math. Phys. 342, (2016),
  491--531.

\bibitem{BD17}
{\sc J.~Breuer and M.~Duits}, {\em {C}entral {L}imit {T}heorems for
  biorthogonal {E}nsembles and asymptotics of recurrence coefficients}.
J. Amer. Math. Soc. 30 No. 1, (2017), 27--66.

\bibitem{CHM16}
{\sc D.~Chafa\"i, A.~Hardy and M.~Ma\"ida}, {\em Concentration for Coulomb gases and Coulomb transport inequalities.}
J. Funct. Anal. 275 (2018), no. 6, 1447--1483.
%\newblock arXiv:1610.00980

\bibitem{CC17}
{\sc C.~Charlier and T.~Claeys},
{\em Thinning and conditioning of the Circular Unitary Ensemble}.
Random Matrices Theory Appl. 6 (2017), no. 2, 1750007, 51 pp. 
%\newblock arXiv:1604.08399

\bibitem{CL95}
{\sc O.~Costin and J.~Lebowitz}, {\em Gaussian fluctuations in random
  matrices}. Phys. Rev. Lett. 75,  (1995), 69--72.

\bibitem{DDMS14}
{\sc D.~S. Dean, P.~L. Doussal, S.~N. Majumdar, and G.~Schehr}, {\em Finite
  temperature free fermions and the {K}ardar-{P}arisi-{Z}hang equation at
  finite time}. Phys. Rev. Lett. 114,  (2015).

\bibitem{Deift99}
{\sc P. Deift},  {\em Orthogonal polynomials and random matrices: a
  {R}iemann-{H}ilbert approach}. Vol.~3 of Courant Lecture Notes in
  Mathematics, Courant Institute of Mathematical Sciences, New York; American Mathematical Society,
  Providence, RI, 1999.

\bibitem{Deift17}
{\sc P. Deift}, {\em Some open problems in random matrix theory and the theory of integrable systems.} II. SIGMA Symmetry Integrability Geom. Methods Appl. 13, (2017), Paper No. 016, 23 pp.



\bibitem{Deift99+}
{\sc P.~Deift, T.~Kriecherbauer, K.T.-R. McLaughlin, S.~Venakides, and
  X.~Zhou}, {\em Uniform asymptotics for polynomials orthogonal with respect to
  varying exponential weights and applications to universality questions in
  random matrix theory}, Comm. Pure Appl. Math. 52 (11),  (1999),
  1335--1425.

\bibitem{DJ16}
{\sc M.~Duits and K.~Johansson}, {\em On mesoscopic equilibrium for linear
  statistics in {D}yson's {B}rownian motion}.
Mem. Amer. Math. Soc. 255 (2018), no. 1222,  118pp
%\newblock arXiv:1312.4295.

\bibitem{Dyson95}
{\sc F.~J. Dyson}, {\em The {C}oulomb fluid and the fifth {P}ainlev{\'e}
  transcendent}, in Chen {N}ing {Y}ang: A Great Physicist of the Twentieth Century Int. Press, Cambridge, MA, 1995,
  131--146.

\bibitem{EK14a}
{\sc L.~Erd\H{o}s and A.~Knowles}, {\em The {A}ltshuler-{S}hklovskii formulas
  for random band matrices {I}: the unimodular case}. Comm. Math. Phy. 333,
  (2015), pp.~1365--1416.

\bibitem{EK14b}
{\sc L.~Erd\H{o}s and A.~Knowles},
 {\em The
  {A}ltshuler-{S}hklovskii formulas for random band matrices {II}: the general
  case}. Ann. Henri Poincar\'{e} 16,  (2015), pp.~709--799.


\bibitem{HK16}
{\sc Y.~He and A.~Knowles}, {\em Mesoscopic eigenvalue statistics of Wigner
  matrices}.
Ann. Appl. Probab. 27 (2017), no. 3, 1510--1550.
%\newblock arXiv:1603.01499, 2016.

\bibitem{HKPV06}
{\sc B.~Hough, M.~Krishnapur, Y.~Peres, and B.~Vir{\'a}g}, {\em Determinantal
  {P}rocesses and {I}ndependence}. Probab. Surv. 3,  (2006), 206--229.

\bibitem{Johansson98}
{\sc K.~Johansson}, {\em On fluctuations of eigenvalues of random {H}ermitian
  matrices}. Duke Math. J. 91,  (1998), 151--204.

\bibitem{Johansson01}
{\sc K.~Johansson}, {\em Universality of the local spacing distribution in
  certain ensembles of {H}ermitian {W}igner matrices}. Comm. Math. Phys. 215,
  (2001), 683--705.


\bibitem{Johansson06}
{\sc K.~Johansson}, {\em Random matrices and determinantal processes}, in
  Mathematical statistical physics, Elsevier B. V., Amsterdam, 2006, 1--55.

\bibitem{Johansson07}
{\sc K.~Johansson}, {\em From {G}umble to {T}racy-{W}idom}. Probab. Theory
  Relat. Fields 138,  (2007), 75--112.


\bibitem{JL15}
{\sc K.~Johansson and G.~Lambert}, {\em Gaussian and non-{G}aussian
  fluctuations for mesoscopic linear statistics in determinantal processes}.
\newblock arXiv:1504.06455.

\bibitem{Konig}
{\sc W.~K\"{o}nig}, {\em Orthogonal polynomial ensembles in probability
  theory}. Probab. Surveys 2,  (2005), 385--447.

\bibitem{Kuijlaars11}
{\sc A.B.J.~Kuijlaars}, {\em Universality} in {\em The Oxford handbook of random matrix theory}.  Oxford Univ. Press, Oxford, 2011, 103--134.


\bibitem{Lambert_a}
{\sc G.~Lambert}, {\em Mesoscopic
  fluctuations for unitary invariant ensembles}.
Electron. J. Probab. 23 (2018), Paper No. 7, 33 pp.
%\newblock arXiv:1510.03641~[v3]

\bibitem{Lambert_b}
{\sc G.~Lambert}, {\em {CLT} for biorthogonal ensembles and related
  combinatorial identities}.
Adv. Math. 329 (2018), 590--648.


\bibitem{LMR_15}
{\sc F.~Lavancier, J.~M{\o}ller, E.~Rubak}, {\em Determinantal point process models and statistical inference}. J. R. Stat. Soc. Ser. B. Stat. Methodol. 77 (2015), no. 4, 853--877. 

\bibitem{LS18}
{\sc T.~Lebl\'e and S. Serfaty}, {\em Fluctuations of two-dimensional Coulomb gases.}
Geom. Funct. Anal. 28 (2018), no. 2, 443--508.
%\newblock arXiv:1609.08088

\bibitem{Pastur06}
{\sc L.A.~Pastur}, {\em Limiting laws of linear eigenvalue statistics for
  {H}ermitian matrix models}. J. Math. Phys. 47,  (2006).

\bibitem{PS11}
{\sc L.A. Pastur and M.~Shcherbina}, {\em Eigenvalue {D}istribution of {L}arge
  {R}andom {M}atrices}. Mathematical Surveys and Monographs 171, Amer. Math.
  Soc., Providence, RI, 2011.

\bibitem{RV07b}
{\sc B.~Rider and B.~Vir{{\'a}}g}, {\em Complex determinantal processes and
  {$H^1$} noise}. Electron. J. Probab. 12, (2007), 1238--1257.

\bibitem{RV07a}
{\sc B.~Rider and B.~Vir{{\'a}}g}. {\em The noise in the
  circular law and the {G}aussian free field}. Int. Math. Res. Not. IMRN
  (2007).

\bibitem{Shcherbina13}
{\sc M.~Shcherbina}, {\em Fluctuations of linear eigenvalue statistics of
  $\beta$-matrix models in the multi-cut regime}. J. Stat. Phys. 151, (2013),
  1004--1034.

\bibitem{Soshnikov}
{\sc A.~Soshnikov}, {\em Determinantal random point fields}. Russian Math.
  Surv.~55,  (2000), 923--975.

\bibitem{Soshnikov00}
{\sc A.~Soshnikov}, 
{\em The {C}entral  {L}imit {T}heorem for local linear statistics in classical compact groups and
  related combinatorial identities}. Ann. Probab.~28,  (2000), 1353--1370.

\bibitem{Soshnikov01}
{\sc A.~Soshnikov}, {\em Gaussian limit for
  determinantal random point fields}. Ann. Probab. 30,  (2001), 1--17.

%\bibitem{SW_13}
%{\sc Sosoe, P., Wong, P.}:
%{\em Regularity conditions in the CLT for linear eigenvalue statistics
   %of Wigner matrices}. Adv. Math. 249, 37--87 (2013)

\bibitem{Spencer11}
{\sc T.~Spencer}, {\em Random banded and sparse matrices} in {\em The Oxford handbook of random matrix theory}.  Oxford Univ. Press, Oxford, 2011, 471--488.


\end{thebibliography}
 \end{document}